\documentclass[11pt]{article}
\usepackage{makeidx}
\usepackage{amsmath,amssymb}
\usepackage{verbatim}

\usepackage{exscale}
\usepackage{tabularx}
\usepackage{amsthm}
\usepackage{latexsym}
\usepackage[usenames,dvipsnames]{color}
\usepackage{graphicx,epsfig}
\usepackage{hyperref}
\usepackage{url}
\usepackage{epstopdf}
\usepackage{tikz}
\usepackage{subfigure}

\usepackage[title,titletoc,toc]{appendix}

\usepackage{natbib}
 \bibpunct[, ]{[}{]}{,}{n}{}{,}%

\newcommand{\mb}[1]{\mbox{\boldmath $#1$}}
\newcommand{\mbs}[1]{{\mbox{\boldmath \scriptsize{$#1$}}}}

\newenvironment{pf}{\noindent {\bf Proof. }}{\hfill $\square$}

\newcommand{\bea}{\begin{eqnarray*}}
\newcommand{\eea}{\end{eqnarray*}}

\newcommand{\vect}[1]{{\boldsymbol #1 }}

\newcommand{\td}[1]{\tilde{#1}}

\newcommand{\card}[1]{\left\lvert #1 \right\rvert}
\newcommand{\bc}{\begin{center}}
\newcommand{\ec}{\end{center}}
\newcommand{\bz}{\vect{z}}

\newcommand{\by}{\vect{y}}

\newcommand{\refs}[1]{$(\ref{#1})$}

\newcommand{\bx}{\vect{x}}

\newcommand{\R}{\mathbb R}

\newcommand{\be}{\begin{equation}}
\newcommand{\ee}{\end{equation}}
\newcommand{\beaa}{\begin{eqnarray*}}
\newcommand{\eeaa}{\end{eqnarray*}}
\newcommand{\ben}{\begin{enumerate}}
\newcommand{\een}{\end{enumerate}}
\newcommand{\db}{\hspace*{\fill}{\zapf o}}

\newcommand{\bpn}{\begin{proposition}\twlsf}
\newcommand{\epn}{\db\end{proposition}}
\newcommand{\bdm}{\begin{displaymath}}
\newcommand{\edm}{\end{displaymath}}
\newcommand{\ba}{\begin{array}}
\newcommand{\ea}{\end{array}}

\newcommand{\st}{\mathop{\rm s.t.}}

\newtheorem{assumption}{Assumption}
\newtheorem{definition}{Definition}
\newtheorem{lemma}{Lemma}
\newtheorem{proposition}{Proposition}
\newtheorem{corollary}{Corollary}
\newtheorem{remark}{Remark}
\newtheorem{theorem}{Theorem}

\newtheorem{example}{Example}
\newcommand{\eps}{\epsilon}

\newcommand{\norm}[1]{\left\lVert#1\right\rVert}

\title{Robust Newsvendor Games with Ambiguity in Demand Distributions}
\author{Xuan Vinh Doan\thanks{DIMAP and ORMS Group, Warwick Business School, University of Warwick, Coventry, CV4 7AL, United Kingdom, xuan.doan@wbs.ac.uk.} \and Tri-Dung Nguyen\thanks{Mathematical Sciences and Business School, University of Southampton, Southampton, SO17 1BJ, United Kingdom, T.D.Nguyen@soton.ac.uk.}}

\date{January 2016}

\begin{document}
\maketitle

\begin{abstract}
We investigate newsvendor games whose payoff function is uncertain due to ambiguity in demand distributions. We discuss the concept of stability under uncertainty and introduce solution concepts for robust cooperative games which could be applied to these newsvendor games. Properties and numerical schemes for finding core solutions of robust newsvendor games are presented.
\end{abstract}
Keywords: Cooperative games; uncertain payoffs; newsvendor games; robust optimization; stability.
%\tableofcontents
%\section{Comments from Dung to Vinh's parts}
%\begin{itemize}
%\item The term `group rational' is often used when describing the stability property ($\sum_{i\in{\cal S}}x_i\geq v({\cal S})$) of the core and hence it might be confusing to use it as an alternative for the term `efficient'. I think just use `efficient' or `efficiency' will be sufficient in this case.
%\item The symbol $\subset$ that we used actually means $\subseteq$ in most of our cases (e.g. in the definition of the core). These two notations are often used interchangeably in many references and it is often OK. However, with the introduction of the formulation (\ref{eq:leastcore}) in our case, it might be better to distinguish these two symbols.
%\item adding \usepackage[dvips]{graphicx,epsfig} somehow make tikz figures not shown properly
%\end{itemize}

\section{Introduction}
\label{sec:intro}
%A1: uncertain payoffs
A joint venture is usually an effective approach for individual players in the market to share costs, reduce risk, and increase the total joint revenue or profit. For example, individual retailers can decide whether to order inventories together and share the profit from selling ordered products later. As a part of the joint venture formation, all players should agree on how to share the joint profit or payoff before the cooperation is established. Cooperative game theory provides a mathematical framework for addressing this problem, which is modeled as \emph{newsvendor centralization games} (or newsvendor games for short). The model of these games has been introduced in Hartman \cite{hartman1994cooperative}. Formally, consider the set ${\cal N}$ of $N$ retailers and let $\td{d}_i\in\mathbb{R}_+$ be the random demand for retailer $i$, $i\in{\cal N}$. In the setting of newsvendor games, we assume that the unit ordering cost $c$ and the unit selling price $p$ are the same for all retailers, $0<c< p$. Given an ordering quantity $y$, the expected profit (or payoff) of retailer $i$ is
\be
\label{eq:optquanti}
v_i(y)=\mathbb{E}_{P_i}\left[p\min\{\td{d}_i,y\}-cy\right],\quad i\in{\cal N}.
\ee
Individual retailer $i$ needs to decide the optimal ordering quantity ${y}^*_i$ to maximize the expected profit (or payoff), $\displaystyle{y}^*_i\in\arg\max_{y\geq 0}v_i(y)$, which is the $(p-c)/p$-quantile of $P_i$, the distribution function of $\td{d}_i$ for all $i\in{\cal N}$. The optimal expected profit is $\displaystyle \bar{v}_i=v_i(y_i^*)=\mathbb{E}_{P_i}\left[p\min\{\td{d}_i,{y}_i^*\}-cy_i^*\right]$.

In the newsvendor games, we are concerned about whether individual retailers should form a coalition to make orders together and share the inventories with each other. For a coalition ${\cal S}\subseteq{\cal N}$, the aggregate demand is $\displaystyle\td{d}({\cal S})=\sum_{i\in{\cal S}}\td{d}_i$ and we assume that the joint distribution $P({\cal S})$ of $\td{d}_i$, $i\in{\cal S}$, is known. Given an ordering quantity $y$, the total expected profit of coalition ${\cal S}$ is
\be
\label{eq:quantbar}
v(y,{\cal S})=\mathbb{E}_{P({\cal S})}\left[p\min\{\td{d}({\cal S}),y\}-cy\right].
\ee
Similarly, coalition ${\cal S}$ of retailers needs to decide the optimal ordering quantity ${y}^*({\cal S})$ to maximize the total expected profit,
\be
\label{eq:optquantbar}
{y}^*({\cal S})\in\arg\max_{y\geq 0}v(y,{\cal S}),
\ee
 which again is the $(p-c)/p$-quantile of the distribution of $\td{d}({\cal S})$. The optimal total expected profit of coalition ${\cal S}$ is
\be
\label{eq:optval}
\bar{v}({\cal S})=v({y}^*({\cal S}),{\cal S})=\mathbb{E}_{P({\cal S})}\left[p\min\{\td{d}({\cal S}),{y}^*({\cal S})\}-cy^*({\cal S})\right].
\ee
In order to guarantee that all retailers are satisfied with the grand coalition $\cal N$, there should be an \emph{allocation} (also called payoff distribution) of the total expected profit $\bar{v}({\cal N})$ among retailers such that no group of retailers has the incentive to form a smaller but better coalition. Newsvendor games concern about the existence of these stable payoff distributions.

Traditionally, the underlying assumption in newsvendor games is that the joint demand distribution is known with certainty. Under this assumption, the expected profits $\bar{v}({\cal S})$ are known for all ${\cal S}\subseteq{\cal N}$. In reality, it can be sometimes difficult to justify this assumption, especially when retailers need to decide whether to cooperate with each other before observing and sharing records of joint demands. If the joint demand distribution is unknown, the expected profits $\bar{v}({\cal S})$ are uncertain. Under the assumption that these payoffs are random variables, several models of stochastic cooperative games have been studied in order to address the uncertainty in coalition payoffs . Charnes and Granot \cite{charnes1973prior} propose a two-stage payoff distribution scheme. In the first stage, some amounts of individual payoffs, which are called prior-payoffs, are promised to the players in the coalition. In the second stage, once the coalition payoff is realized, there could be adjustments in individual payoffs to guarantee that the payoff distribution is feasible and that any objections among the players are minimized. This approach assumes risk-neutral behaviors among players. Suijs et al. \cite{suijs1999cooperative} study a different payoff distribution scheme for stochastic cooperative games using preference orders for random payoffs, which can handle different types of risk behavior. An individual payoff under this distribution scheme consists of two parts. The first part is a monetary exchange between players and the second part is a fraction of the random coalition payoff, which can depend on the action taken by the coalition. This distribution scheme works well for some applications such as the insurance game discussed in Suijs et al. \cite{suijs1998stochastic}. More recently, Uhan \cite{uhan2013stochastic} generalizes this payoff distribution scheme for stochastic linear programming games with applications in inventory centralization and network fortification. Timmer et al. \cite{timmer2005convexity} argue that risk-covering monetary compensation is not needed in other applications. It is indeed a reasonable assumption that players (e.g., retailers or firms) in the examples discussed therein do not feel the need to pay upfront in order to benefit from a greater share of the profit in the future. Timmer et al. \cite{timmer2005convexity} then propose a similar payoff distribution scheme for stochastic cooperative games without monetary exchange in which a multiple of the random coalition payoff is allocated to each individual. Fern\'andez et al. \cite{fernandez2002cores} also consider this payoff distribution scheme as a special case while investigating general stochastic payoff distributions using stochastic orders for stochastic cooperative games.

%A3: other models with uncertain payoffs
In addition to these models of stochastic cooperative games, there are different models that address different aspects of cooperative games under uncertainty. Lehrer \cite{lehrer2003allocation} investigates payoff distribution processes in repeated deterministic cooperative games. The payoff allocation at any time depends on past allocations and current coalition payoffs. Lehrer \cite{lehrer2003allocation} shows that these payoff distribution processes converge to some well-known solutions of cooperative games under appropriate allocation rules. Dror et al. \cite{dror2008dynamic} extend these allocation rules and payoff distribution processes for repeated stochastic cooperative games and apply them to dynamic newsvendor realization games. Bauso and Timmer \cite{bauso2009robust}, on the other hand, study dynamic cooperative games under the setting of a family of games whose coalition values are uncertain. Toriello and Uhan \cite{toriello2015dynamic} investigate the dynamic linear programming games with risk-averse players. The study of cooperative games with incomplete information addresses other aspects of uncertainty in cooperative games. Ieong and Shoham \cite{ieong2008bayesian} develop solution concepts of Bayesian cooperative games such as \emph{ex ante}, \emph{ex interim}, and \emph{ex post} cores under the assumption that all players have a common prior over the set of possible states of the world but that each of them has some private information or belief about the true state of the world. More recently, Li and Conitzer \cite{li2014cooperative} focus on solution concepts that maximize the probability of \emph{ex post} stability assuming all players know the distribution of the random state of the world. Clearly, in this setting, private information of players does not play an important role since only \emph{ex post} stability is considered. Chalkiadadis and Boutilier \cite{chalkiadakis2004bayesian} study Bayesian coalition formation problems when players have different beliefs about contributions of other players to coalition payoffs, which makes coalition payoffs uncertain. Forges and Serrano \cite{forges2013cooperative} discuss different models of cooperative games with incomplete information that handle externality; i.e., the effects of actions and information of players outside a coalition on the coalition payoff and payoffs of individual players in the coalition.

%A4: our approach
In this paper, we focus on newsvendor games with uncertain expected profits (or payoffs) under a similar setting to those of stochastic cooperative games discussed above; that is, uncertainty is captured and represented only in the payoff functions of the games. Our research, however, is different from existing models of stochastic cooperative games given the fact it is difficult to obtain probabilistic information on uncertain expected profits $\bar{v}({\cal S})$, ${\cal S}\subseteq{\cal N}$. This assumption is particularly useful for newsvendor games considered in this paper in which the uncertainty of the payoffs is modeled through an uncertainty set derived from the ambiguity of the joint demand distribution.
\subsection*{Contributions and paper outline}
In this paper, we focus on the newsvendor games with ambiguity in the joint demand distribution. Specifically, our contributions and the structure of the paper are as follows:
\ben
\item[(1)] Section~\ref{ssec:uncertain} develops a framework of cooperative games with uncertain payoffs which could be applied to newsvendor games under distributional ambiguity. We call these games \emph{robust cooperative games}. Solution concepts such as imputation, core, and game balancedness of robust cooperative games are defined and discussed.
\item[(2)] Section~\ref{stoc.newsvendor} studies newsvendor games with ambiguity in demand distributions using the framework of robust cooperative games. We focus on the existence of rational and stable payoff distributions of these \emph{robust newsvendor games}. We also discuss the computational aspect of finding core solutions for robust newsvendor games and provide some numerical results. %A demonstration of the properties of the robust solutions and computational results are provided in Section~\ref{sec:numerical_results}. We present a constraint generation method for speeding up the computation of the robust payoff distribution in larger games with more than $30$ retailers.
\een
%Section~\ref{sec:conclusion} concludes and provides a discussion of future directions.
\section{Deterministic Cooperative Games}
\label{sec:coop}
Before presenting the model of cooperative games with uncertain payoffs that could be applied for newsvendor games under distributional ambiguity, we introduce some solution concepts of deterministic cooperative games whose payoffs are known with certainty. Consider a set of $N$ players, ${\cal N}=\{1,\ldots,N\}$, and a function $v:2^{\cal N}\rightarrow\R$, which is called the \emph{characteristic function}, such that $v(\emptyset)=0$. We call ${\cal G}=({\cal N},v)$ a \emph{cooperative game}. Utility (payoff) is assumed to be transferable, i.e., for any coalition ${\cal S}\subseteq{\cal N}$, its total payoff is completely defined as $v({\cal S})$, which can be transferred freely among its members. Under the assumption that the joint demand distribution is known with certainty, the newsvendor games previously introduced is indeed a deterministic cooperative game, $\bar{\cal G}=({\cal N},\bar{v})$, where $\bar{v}$ is defined as in \refs{eq:optval}.

Given a cooperative game $({\cal N},v)$, we are interested in finding an allocation $\bx\in\R^N$ to distribute the total payoff $v({\cal N})$ among individual players. An allocation (also called a payoff distribution) $\bx$ is \emph{efficient} if
\be
\label{eq:grouprational}
\sum_{i\in {\cal N}}x_i=v({\cal N}).
\ee
An important question regarding cooperative games is whether players are willing to join the grand coalition $\cal N$. Necessarily, there should be an allocation $\bx\in\R^N$ with which each individual player is better off as compared to his/her standalone payoff. An allocation $\bx$ is called \emph{individually rational} if
\be
\label{eq:indrational}
x_i\geq v(\{i\}),\quad\forall\,i\in{\cal N}.
\ee
\begin{definition}
\label{def:imputation}
An \emph{imputation} is an allocation that is both efficient and individually rational. The set of imputations of a cooperative game $({\cal N},v)$ is written as follows:
\be
\label{eq:impu}
\emph{impu}({\cal N},v):=\left\{\bx\in\R^N\,\mid\,\sum_{i\in {\cal N}}x_i=v({\cal N}),\,x_i\geq v(\{i\}),\,\forall\,i\in{\cal N}\right\}.
\ee
\end{definition}
If the characteristic function $v$ is \emph{super-additive}, i.e., for any two disjoint coalitions ${\cal S}_1$ and ${\cal S}_2$, ${\cal S}_1\cap{\cal S}_2=\emptyset$, $v({\cal S}_1)+v({\cal S}_2)\leq v({\cal S}_1\cup{\cal S}_2)$, it is clear that there always exists at least one imputation or $\mbox{impu}({\cal N},v)\neq\emptyset$.

Individual rationality is not sufficient to guarantee that some players would prefer the grand coalition $\cal N$ to a smaller coalition ${\cal S}\subsetneq{\cal N}$. An allocation is called \emph{stable} with respect to a coalition $\cal S$ if $\displaystyle\sum_{i\in{\cal S}}x_i\geq v({\cal S})$.
\begin{definition}
\label{def:core}
The set of efficient allocations that are stable with respect to all coalitions ${\cal S}\subseteq{\cal N}$ is called the \emph{core},
\be
\label{eq:core}
\emph{core}({\cal N},v):=\left\{\bx\in\R^N\,\mid\,\sum_{i\in {\cal N}}x_i=v({\cal N}),\,\sum_{i\in{\cal S}}x_i\geq v({\cal S}),\,\forall\,{\cal S}\subsetneq{\cal N}\right\}.
\ee
\end{definition}
It is obvious that $\mbox{core}({\cal N},v)\subseteq\mbox{impu}({\cal N},v)$. However, the core might not exist. There are different solution concepts that take this issue into account. Given a parameter $\eps\geq 0$, the $\eps$-\emph{core} is defined as follows:
\be
\label{eq:ecore}
\eps\mbox{-core}({\cal N},v):=\left\{\bx\in\R^N\,\mid\,\sum_{i\in {\cal N}}x_i=v({\cal N}),\,\sum_{i\in{\cal S}}x_i\geq v({\cal S})-\eps,\,\forall\,{\cal S}\subsetneq{\cal N}\right\}.
\ee
It is clear that $\eps$-core is non-empty for a given large enough $\eps$. The \emph{least core} is the non-empty $\eps$-core with the smallest value of $\eps$, which is called the \emph{least core value}. The least core value $\eps({\cal N},v)$ is the optimal value of the following linear program:
\be
\label{eq:leastcore}
\ba{rll}
\displaystyle\eps({\cal N},v)=\min_{\vect{x},\eps\geq 0} & \eps\\
\st & \displaystyle\sum_{i\in{\cal N}}x_i=v({\cal N}),\\
& \displaystyle\sum_{i\in{\cal S}}x_i\geq v({\cal S})-\eps, & \forall\,{\cal S}\subseteq{\cal N}.
\ea
\ee
%If we are able to solve \refs{eq:leastcore}, the non-emptiness of the core can be determined and if the core is non-empty, an allocation in the core is also determined; otherwise, we achieve an allocation in the least core.
If we are able to solve Problem~\refs{eq:leastcore} and obtain $\eps({\cal N},v) = 0$, then the core is non-empty. Schulz and Uhan \cite{Schulz13} consider Problem \refs{eq:leastcore} in a slightly different way:
\be
\label{eq:stability}
\ba{rll}
\displaystyle s({\cal N},v)=\min_{\vect{x},\eps} & \eps\\
\st & \displaystyle\sum_{i\in{\cal N}}x_i=v({\cal N}),\\
& \displaystyle\sum_{i\in{\cal S}}x_i\geq v({\cal S})-\eps, & \forall\,{\cal S}\subsetneq{\cal N}, {\cal S} \neq \emptyset.
\ea
\ee
In general, we have $\eps({\cal N},v)=\max\{s({\cal N},v),0\}$, i.e., $s({\cal N},v)$ coincides with the least core value if it is non-negative. If $s({\cal N},v)<0$, the absolute value $\card{s({\cal N},v)}$ can be interpreted as the maximum possible increase in all coalition payoffs $v({\cal S})$, ${\cal S}\subsetneq{\cal N}$, under the condition that at least one efficient and stable allocation still exists. %We define $\sigma({\cal N},v)=\max\{-s({\cal N},v),0\}$ to be the \emph{stability value} of the game $({\cal N},v)$. Intuitively, the larger the stability value is, the more stable the grand coalition is under changes in the payoffs of smaller coalitions. The stability value can be computed given a fixed characteristic function and proves to be a useful concept when we later discuss cooperative games with uncertain characteristic functions. We have so far introduced the concepts of cooperative games with a fixed characteristic function $v$.

For the newsvendor game $({\cal N},\bar{v})$, the characteristic function $\bar{v}$ is super-additive. Indeed, we have
$$
\bar{v}({\cal S})=\mathbb{E}_{P({\cal S})}\left[p\min\{\td{d}_i,{y}^*({\cal S})\}-cy^*({\cal S})\right]=(p-c){y}^*({\cal S})-p\,\mathbb{E}_{P({\cal S})}\left[\left({y}^*({\cal S})-\td{{d}}({\cal S})\right)^+\right],
$$
where $x^+=\max\{x,0\}$. Consider two disjoint coalitions ${\cal S}_1$ and ${\cal S}_2$, ${\cal S}_1\cap{\cal S}_2=\emptyset$; we have $\td{d}({\cal S}_1\cup{\cal S}_2)=\td{d}({\cal S}_1)+\td{d}({\cal S}_2)$. Using \refs{eq:optquantbar} and the fact that $(x+y)^+\leq x^++y^+$, we have
$$
\ba{rl}
\bar{v}({\cal S}_1\cup{\cal S}_2)&\geq\,(p-c)({y}^*({\cal S}_1)+{y}^*({\cal S}_2))-p\,\mathbb{E}_{P({\cal S}_1\cup{\cal S}_2)}\left[\left({y}^*({\cal S}_1)+{y}^*({\cal S}_2)-\td{{d}}({\cal S}_1\cup{\cal S}_2)\right)^+\right]\\
&\geq\,(p-c){y}^*({\cal S}_1)-p\,\mathbb{E}_{P({\cal S}_1)}\left[\left({y}^*({\cal S}_1)-\td{{d}}({\cal S}_1)\right)^+\right]+\vspace{5pt}\\
&\phantom{\geq\,\,}\,(p-c){y}^*({\cal S}_2)-p\,\mathbb{E}_{P({\cal S}_2)}\left[\left({y}^*({\cal S}_2)-\td{{d}}({\cal S}_2)\right)^+\right]\\
&=\,\bar{v}({\cal S}_1)+\bar{v}({\cal S}_2).
\ea
$$
This shows that newsvendor games always have imputations. For normally distributed demands, Hartman et al. \cite{hartman2000cores} show that the cores of these games are non-empty. \citet{chen2009stochastic} use stochastic linear programming duality to show the non-emptiness of the core in inventory centralization games and to provide a computational method for finding one. M\"uller et al. \cite{mueller02} prove that every newsvendor game has a non-empty core, i.e., for all possible distributions of random demands. Montrucchio and Scarsini \cite{montrucchio07} show how to construct an allocation in the core of newsvendor games. In the next section, we develop a model of cooperative games with uncertain characteristic functions that could be used to study newsvendor games with the ambiguity in demand distributions.

\section{Robust Cooperative Games}
\label{ssec:uncertain}
When the payoffs are uncertain and assumed to be random variables (with finite expectation), the model of stochastic cooperative games proposed by Suijs et al. \cite{suijs1999cooperative} can be applied. Suijs et al. \cite{suijs1999cooperative} define a stochastic cooperative game as a tuple $({\cal N},{\cal A},\tilde{v},\succsim)$, where ${\cal N}$ is the set of players, ${\cal A}({\cal S})$ is the set of possible actions that coalition $\cal S$, ${\cal S}\subseteq{\cal N}$, can take. $\tilde{v}(a,{\cal S})$ is the random payoff obtained by coalition $\cal S$ if action $a\in{\cal A}({\cal S})$ is taken. The inclusion of action sets ${\cal A}({\cal S})$ in this model of stochastic cooperative games reflects the fact that players can take different actions when facing uncertainty, which is different in the setting of deterministic cooperative games in which everything, including actions of players, is known with certainty. For newsvendor games, if retailers in a coalition ${\cal S}$ do not know their joint demand distribution, there is no ordering quantity to achieve the ``optimal'' expected profit and they can choose to pick an ordering quantity $y$ from a set of possible ordering quantities ${\cal Y}({\cal S})$. Clearly, the expected profit of a coalition ${\cal S}$, $\tilde{v}(y,{\cal S})$, now depends on the ordering quantity $y$. Finally, $\succsim_i$ describes the preference of player $i$, $i\in{\cal N}$, over different random variables. For two random payoff $\tilde{v}$ and $\tilde{v}'$, $\tilde{v}\succsim_i\tilde{v}'$ means that player $i$ prefers $\tilde{v}$ over $\tilde{v}'$. This order relation of stochastic payoffs is needed for this model of stochastic cooperative games to replace the straightforward comparison among deterministic payoffs in the deterministic setting. Suijs et al. \cite{suijs1999cooperative} discuss several preference relations such as quantile-based relations and stochastic dominance for which the existence of core allocations is studied in detail.

For newsvendor games, the assumption of a known joint demand distribution in the deterministic setting can be considered to be rather strong given the fact that individual retailers usually collect historical demands independently before they join any coalition. In order to make the problem more realistic, we can assume that some (multivariate) marginal distributions are known instead of the joint distribution. For example, it is more reasonable to assume the knowledge of the joint demand distribution of a subset of retailers that are located close to each other and hence likely to serve customers from the same area. Under this assumption of distributional ambiguity, the joint demand distribution is unknown and belongs to an ambiguity set. For a coalition $\cal S$, given an ordering quantity $y$, the total expected profit $\td{v}({y,\cal S})=\mathbb{E}_{\td{P}({\cal S})}\left[p\min\{\td{d}_i,{y}\}-cy({\cal S})\right]$ is uncertain and belongs to an uncertainty set if the (marginal) distribution $\td{P}({\cal S})$ is not known with certainty. In general, the characteristic function $\td{v}$ of newsvendor games is uncertain if there is ambiguity in demand distributions. In order to apply the model of stochastic cooperative games discussed above for these newsvendor games, we would need to assume that probability distributions of random payoffs $\td{v}(y,{\cal S})$ are known. However, it is difficult to impose this probabilistic assumption, given that the randomness of $\td{v}(y,{\cal S})$ comes from the uncertainty of the probability distribution $\td{P}({\cal S})$. In order to address this issue, we propose a new model of general cooperative games with uncertain characteristic functions, which is different from the model of stochastic cooperative games proposed by Suijs et al. \cite{suijs1999cooperative}.

Given a coalition ${\cal S}\subseteq{\cal N}$ and an action $a\in{\cal A}({\cal S})$, instead of considering $\td{v}(a,{\cal S})$ as a random variable, we consider it as an \emph{uncertain} payoff $v_u(a,{\cal S})$, where $u\in{\cal U}$ and ${\cal U}$ is an uncertainty set. %The set of all characteristic functions is denoted by ${\cal V}({\cal U})=\{v_u\,:\,{\cal A}\times 2^{{\cal N}}\rightarrow\mathbb{R}\}$ and let ${\cal V}_{{\cal U}}(a,{\cal S})=\{v_u(a,{\cal S})\,:\,u\in{\cal U}\}$ for ${\cal S}\subseteq{\cal N}$ and $a\in{\cal A}({\cal S})$.
Let us denote ${\cal V}_{{\cal U}}(a,{\cal S})=\{v_u(a,{\cal S})\,:\,u\in{\cal U}\}$ for ${\cal S}\subseteq{\cal N}$ and $a\in{\cal A}({\cal S})$. We also denote ${\cal V}({\cal U})=\{v_u\,:\,{\cal A}\times 2^{{\cal N}}\rightarrow\mathbb{R},u\in{\cal U}\}$ as the set of all characteristic functions. This setting is appropriate for the newsvendor games with the ambiguity in demand distribution in which the joint demand distribution $P$ belongs to an ambiguity set ${\cal P}$. We define an instance of cooperative games with uncertain characteristic functions by the tuple $({\cal N},{\cal A},{\cal V}({\cal U}))$. In order to define the main concepts of these cooperative games $({\cal N},{\cal A},{\cal V}({\cal U}))$, we need first to select a payoff distribution scheme and characterize when sub-coalitions would break away from the grand coalition.

There are many different ways how payoffs can be allocated under uncertainty. When payoffs are assumed to be random variables, the most general payoff distribution scheme with random allocation for stochastic cooperative games is studied by Fern\'andez et al. \cite{fernandez2002cores}. For a coalition ${\cal S}$, they consider efficient allocations as random vectors $\td{\bx}$ that satisfy
\be
\label{eq:sefficiency}
\sum_{i\in{\cal S}}\td{x}_i=\td{v}(a,{\cal S}),
\ee
assuming $a\in{\cal A}({\cal S})$ is the chosen action for the coalition. Suijs et al. \cite{suijs1999cooperative} argue that this general stochastic payoff distribution scheme generates a very large class of allocations and results in computational difficulties. They propose a more restrictive payoff distribution scheme that can be represented by a pair $(\mb{d},\mb{r})\in\R^{\card{\cal S}}\times\R^{\card{\cal S}}$ such that $\displaystyle\sum_{i\in{\cal S}}d_i=\mathbb{E}\left[\tilde{v}(a,{\cal S})\right]$ and $\displaystyle\sum_{i\in{\cal S}}r_i=1$, $r_i\geq 0$ for all $i\in{\cal S}$. The stochastic allocation is defined as $\td{x}_i=d_i+r_i\left(\td{v}(a,{\cal S})-\mathbb{E}\left[\tilde{v}(a,{\cal S})\right]\right)$ for all $i\in{\cal S}$. Note that this allocation scheme require probabilistic information, i.e., expectations of random characteristic functions. When the characteristic functions are deterministic, this allocation scheme coincides with the classical allocation of deterministic cooperative games, $x_i=d_i$ for all $i\in{\cal S}$, where $\displaystyle\sum_{i\in{\cal S}}d_i=\mathbb{E}\left[\tilde{v}(a,{\cal S})\right]=v({\cal S})$. Timmer et al. \cite{timmer2005convexity} propose another payoff distribution scheme that is represented by a vector $\mb{z}\in\R^S$ of \emph{multiples}. The stochastic allocation under this scheme is $\td{x}_i=\td{v}(a,{\cal S})\cdot z_i$ for all $i\in{\cal S}$. Here, it is not required for $z_i$ to be non-negative for $i\in{\cal S}$. In order to guarantee that the stochastic allocation is efficient, it is required that $\displaystyle\sum_{i\in{\cal S}}z_i=1$. Compared to the payoff distribution scheme in Suijs et al. \cite{suijs1999cooperative}, Timmer et al. \cite{timmer2005convexity} argue that their proposed scheme does not allow monetary compensations to cover risk, which can be considered as a reasonable setting for several applications such as newsvendor games discussed above. Note that this scheme can also be seen as a direct extension of that of deterministic cooperative games since any allocation of deterministic cooperative games can be written in terms of multiples of $v({\cal S})$, assuming $v({\cal S})\neq 0$. For the proposed model of cooperative games with uncertain characteristic functions, we face the same issues when defining payoff distribution schemes due to the uncertainty of the characteristic functions. Instead of having a random payoff $\tilde{x}_i$ as in stochastic cooperative games, the uncertain payoff of player $i$ is now represented by the set ${\cal X}_i({\cal U})=\{x_i(u)\,:\,u\in{\cal U}\}\subset\mathbb{R}$. For a coalition ${\cal S}$, we obtain the set of payoff vectors $\mb{{\cal X}}_{{\cal S}}({\cal U})=\{\bx_{{\cal S}}(u)\,:\,u\in{\cal U}\}\subset\mathbb{R}^{\card{{\cal S}}}$, where $\bx_{{\cal S}}(u)=\{x_i(u)\,:\,i\in{\cal S}\}$. Given that our main application is the newsvendor game in which risk-covering monetary compensation is not very important, we adopt a similar payoff distribution scheme to the one proposed by Timmer et al. \cite{timmer2005convexity} for our model of cooperative games with uncertain characteristic functions. Formally, for a coalition ${\cal S}$, a payoff distribution scheme is represented by $\bz \in\R^{\card{\cal S}}$ such that for a given action $a\in{\cal S}$, the uncertain allocation for each player $i$, $i\in{\cal S}$ is
\be
\label{eq:scheme}
x_i(u)=v_u(a,{\cal S})\cdot z_i,\quad\forall\,u\in{\cal U}.
\ee
A payoff distribution scheme $\bz$ is \emph{efficient} if $\displaystyle \sum_{i\in{\cal S}}z_i=1$.
%\be
%\label{eq:efficient}
%\sum_{i\in{\cal S}}z_i=1.
%\ee
%This modified scheme is a proper extension of that of deterministic cooperative games for any value of the characteristic function.
This concept of efficiency indicates that the uncertain allocation $\bx$ is efficient for all realizations of the uncertain characteristic function $v_u\in{\cal V}({\cal U})$.

Now, similar to stochastic cooperative games, the concepts of rationality and stability rely on both the payoff distribution scheme and the action taken by the grand coalition, the \emph{decision} $(a,\bz)$. They also depend on how we characterize the preference of individuals and sub-coalitions with respect to staying in or breaking away from the grand coalition. In the model of stochastic cooperative games proposed by Suijs et al. \cite{suijs1999cooperative}, order relations of stochastic payoffs are used to define these stay-in/break-away preferences. Under the setting of uncertain payoffs, there is no distributional information and we need to define these preferences based on the realization of uncertain parameters. The resulting uncertain allocation for each player $i$, $i\in{\cal N}$ is $x_i(u)=v_u(a,{\cal N})\cdot z_i$, $u\in{\cal U}$ for all $i\in{\cal N}$. If he/she acts alone with the action $\alpha\in{\cal A}(\{i\})$, the uncertain payoff is $v_u(\alpha,\{i\})$. If he/she joins a smaller sub-coalition ${\cal S}\subsetneq{\cal N}$ whose decision is $(\hat{a},\hat{\mb{z}})$, the uncertain allocation would be $v_u(\hat{a},{\cal S})\cdot\hat{z}_i$. Given a realization of $u\in{\cal U}$, each player $i$ can compare these different realized payoffs directly and the sub-coalition would have the incentive to break away if all of its members are better off as in the deterministic setting. By looking at the break away incentive given a realization of uncertain parameters, we can define stay-in/break-away preferences of sub-coalitions (and individual players) under the uncertain setting as follows.
\begin{definition}
\label{def:preference}
A sub-coalition ${\cal S}$ has the incentive to break away from the grand coalition if there exists a decision $(\hat{a},\hat{\mb{z}})$ which guarantees that \emph{all} players in the sub-coalition are better off by leaving given \emph{at least} one realization $u\in{\cal U}$, i.e., $v_u(\hat{a},{\cal S})\cdot\hat{z}_i>x_i(u)$ for all $i\in{\cal S}$.
\end{definition}
In other words, a sub-coalition has the incentive to stay in the grand coalition if no matter what decision it takes, there will be at least one player in the sub-coalition is better off by staying in the grand coalition in any realization of uncertain parameters. Applying this definition for individual players, we can say that a player has the incentive to break away if there exists an action $\alpha\in{\cal A}(\{i\})$  which guarantees that he/she will be better off by leaving the grand coalition for at least one realization of $u\in{\cal U}$, i.e., $v_u(\alpha,\{i\})>x_i(u)$. Note that the idea of defining these realization-based concepts is similar to those in Bitran \cite{bitran1980linear}, which are used to define necessary and sufficient solution concepts for multi-objective optimization. We are now ready to define main solution concepts for cooperative games with uncertain characteristic functions $({\cal N},{\cal A},{\cal V}({\cal U}))$ using the proposed payoff distribution scheme and the above definition of break away incentive.

A decision $(a,\bz)$ is \emph{individually rational} if there is no individual player who has the incentive to break away. We can then define the concept of imputations.
\begin{definition}
\label{def:rimputation}%%%%%
An imputation of the cooperative game $({\cal N},{\cal A},{\cal V}({\cal U}))$ is an individually rational decision $(a,\bz)$ whose payoff distribution scheme $\bz$ is efficient.
\end{definition}
Similar to the deterministic setting, let $\mbox{impu}({\cal N},{\cal A},{\cal V}({\cal U}))$ denote the set of all imputations of the cooperative game $({\cal N},{\cal A},{\cal V}({\cal U}))$.

We now define the concept of stability. A decision $(a,\bz)$ is \emph{stable} if there is no sub-coalition ${\cal S}\subsetneq{\cal N}$ that has the incentive to break away and we have the following definition of the cores of cooperative games $({\cal N},{\cal A},{\cal V}({\cal U}))$.
\begin{definition}
\label{def:score}
The core of the cooperative game $({\cal N},{\cal A},{\cal V}({\cal U}))$ is the set of all stable decisions with efficient payoff distribution schemes and is denoted by $\emph{core}({\cal N},{\cal A},{\cal V}({\cal U}))$.
\end{definition}

The above definition of stability indicates that all sub-coalitions of a stable grand coalition necessarily do not have the incentive to break away no matter how $u\in {\cal U}$ is realized.  
%Dung comment: It can be shown that this statement is equivalent to the definition 3 on stability.
This follows the ``immunized-against-uncertainty'' principle of robust optimization (see Ben-Tal et al. \cite{ben2009robust} and references therein); therefore, we call cooperative games $({\cal N},{\cal A},{\cal V}({\cal U}))$ with the above definition of break away incentive \emph{robust cooperative games}.

We are now ready to characterize the existence of imputations and core decisions, i.e., decisions belong to the core, of robust cooperative games $({\cal N},{\cal A},{\cal V}({\cal U}))$. We start by making the following assumptions.

\begin{assumption}
\label{as:pos}
\indent
\begin{itemize}
\item[(i)] For any player $i\in{\cal N}$, there exists an action such that his/her payoff is always non-negative, i.e., there exists $a \in {\cal A}(\{i\})$ such that $v_u(a,\{i\})\geq 0$ for all $u\in{\cal U}$.
\item[(ii)] The payoff of the grand coalition $\cal N$ is always positive, i.e., $v_u(a,{\cal N})>0$ for all $a\in{\cal A}({\cal N})$ and $u\in{\cal U}$.
\end{itemize}
\end{assumption}
Similar to Timmer et al. \cite{timmer2005convexity}, these assumptions emphasize the fact that we are focusing on profit games with possible nonnegative payoff for each individual sub-coalition and it is indeed worth considering the grand coalition given that its profit is always positive. In other words, Assumption~\ref{as:pos}(ii) implies that we should only consider the set of actions ${\cal A}({\cal N})$ which create positive profits under any circumstances for the grand coalition. For robust newsvendor games that we are going to discuss in Section~\ref{stoc.newsvendor}, these assumptions are easily satisfied in general.

%Assumption~\ref{as:pos}(i) automatically holds if we assume that, for any coalition, there exists a `null' action such that the coalition always receives zero profit. In the robust newsvendor games to be described in Section~\ref{stoc.newsvendor}, this null action is essentially the action of ordering nothing, i.e. do not do the business, and hence always receiving a payoff of zero.

%\begin{assumption}
%\label{as:null:action}
%For any coalition $\cal S$, there exists an action such that the payoff of the coalition is always non-negative, i.e., there exists $\alpha \in {\cal A}({\cal S})$ such that $v_u(\alpha,{\cal S})\geq 0$ for all $u\in{\cal U}$.
%\end{assumption}

We now state the following result on existence conditions for imputations of robust cooperative games.
\begin{theorem}
\label{prop:impucond}
Given a robust cooperative game $({\cal N},{\cal A},{\cal V}({\cal U}))$, an imputation exists if and only if there exists an action $a\in{\cal A}({\cal N})$ such that
\be
\label{eq:impucond}
\sum_{i\in{\cal N}} v_{\max}(a,\{i\}) \leq 1,
\ee
where $\displaystyle v_{\max}(a,\{i\}) = \max_{u\in{\cal U}}\left\{\frac{\displaystyle\max_{\alpha_i\in{\cal A}(\{i\})}v_u(\alpha_i,\{i\})}{v_u(a,{\cal N})}\right\}$.
\end{theorem}
In order to prove the theorem, we need the following lemma.
\begin{lemma}
\label{prop:individual:rationality:cond}
Given a robust cooperative game $({\cal N},{\cal A},{\cal V}({\cal U}))$ and a decision $(a,\bz)$, a player $ i \in {\cal N}$ has no incentive to break away if and only if
$\displaystyle z_i \geq v_{\max}(a,\{i\})$.
\end{lemma}

\begin{pf}
We have: $x_i(u)=v_u(a,{\cal N})\cdot z_i$ for all $i\in{\cal N}$ and $u\in{\cal U}$. Player $i\in{\cal N}$ does not have the incentive to break away if and only if %${\cal X}_i({\cal U})\succsim_n{\cal V}_{\cal U}(\alpha_i,\{i\})$ for all $\alpha_i\in{\cal A}(\{i\})$, i.e., for all $u\in{\cal U}$,
$$
v_u(\alpha_i,\{i\})\leq v_u(a,{\cal N})\cdot z_i,\quad\,\forall\,\alpha_i\in{\cal A}(\{i\}),~u\in{\cal U}.
$$
Under Assumption \ref{as:pos}(ii), this holds if and only if $\displaystyle z_i\geq\max_{\alpha_i\in{\cal A}(\{i\})}\frac{v_u(\alpha_i,\{i\})}{v_u(a,{\cal N})}$ for all $u\in{\cal U}$, i.e.,
$$
z_i\geq \max_{u\in{\cal U}}\left\{\frac{\displaystyle\max_{\alpha_i\in{\cal A}(\{i\})}v_u(\alpha_i,\{i\})}{v_u(a,{\cal N})}\right\}=v_{\max}(a,\{i\}).
$$
\end{pf}

We are now ready to prove Theorem \ref{prop:impucond}.

\begin{pf}
Suppose an imputation $(a,\bz)$ exists. According to the definition of imputations, there is no player $i\in{\cal N}$ who has the incentive to break away. Thus, according to Lemma~\ref{prop:individual:rationality:cond}, we have:
$$
z_i\geq v_{\max}(a,\{i\}),\quad\,\forall\,i\in{\cal N}.
$$
Now, $(a,\bz)$ is a decision of the grand coalition; therefore, $\displaystyle\sum_{i\in{\cal N}}z_i=1$. Summing over all $i\in{\cal N}$ the above inequality, we then achieve condition \refs{eq:impucond}.

Now, suppose condition \refs{eq:impucond} holds. Let $\displaystyle\eps = 1-\sum_{i\in{\cal N}} v_{\max}(a,\{i\}) \geq 0$ and define
$$
z_i=v_{\max}(a,\{i\})+\frac{\eps}{N},\quad\,\forall\,i\in{\cal N}.
$$
We will show that $(a,\bz)$ is an imputation. Clearly, $\bz$ is efficient, i.e., $\displaystyle\sum_{i\in{\cal N}}z_i=1$ given the definition of $z_i$ and $\eps$. Now we have: for all $u\in{\cal U}$, $\displaystyle z_i\geq v_{\max}(a,\{i\}) \geq \frac{\displaystyle\max_{\alpha_i\in{\cal A}(\{i\})}v_u(\alpha_i,\{i\})}{v_u(a,{\cal N})}$ for all $i\in{\cal N}$ since $\eps\geq 0$. Thus, under Assumption \ref{as:pos}(ii),
$$
v_u(a,{\cal N})\cdot z_i\geq \max_{\alpha_i\in{\cal A}(\{i\})}v_u(\alpha_i,\{i\}),\quad\,\forall\,i\in{\cal N},\,u\in{\cal U}.
$$
It shows that for all $u\in{\cal U}$,
$$
x_i(u)\geq v_u(\alpha_i,\{i\}),\quad \forall\,i\in{\cal N},\,\alpha_i\in{\cal A}(\{i\}).
$$
%or equivalently, ${\cal X}_i({\cal U})\succsim_n{\cal V}_{\cal U}(\alpha_i,\{i\})$ for all $\alpha_i\in{\cal A}(\{i\})$. 
Thus, there is no player $i$ who has the incentive to break away; that is $(a,\bz)$ is individually rational, which implies that $(a,\bz)$ is an imputation.
\end{pf}

Similar to the deterministic and stochastic cooperative games, the existence of core decisions is related to the concept of balancedness. A map $\mu:2^{{\cal N}}\setminus\emptyset\rightarrow[0,+\infty)$ is called a \emph{balanced} map if $\displaystyle\sum_{{\cal S}\subsetneq{\cal N}}\mu({\cal S})\cdot\mb{e}_{{\cal S}}=\mb{e}_{{\cal N}}$, where, for all ${\cal S}\subseteq{\cal N}$, $\mb{e}_{{\cal S}}\in\{0,1\}^N$ with $\left(e_{{\cal S}}\right)_i=1$ if and only if $i\in{\cal S}$. We are now ready to define the balanced robust cooperative game.
\begin{definition}
A robust cooperative game $({\cal N},{\cal A},{\cal V}({\cal U}))$ is called balanced if there exists an action $a\in{\cal A}({\cal N})$ such that for all balanced map $\mu$,
\be
\label{eq:balanced}
\sum_{{\cal S}\subsetneq{\cal N}}\mu({\cal S})v_{\max}(a,{\cal S})\leq 1,
\ee
where $\displaystyle v_{\max}(a,{\cal S})=\max_{u\in{\cal U}}\left\{\frac{\displaystyle\max_{\alpha_{{\cal S}}\in{\cal A}({\cal S})}v_u(\alpha_{{\cal S}},{\cal S})}{v_u(a,{\cal N})}\right\}$.
\end{definition}
This definition of balanced robust cooperative games matches the definition of balanced deterministic games when $\card{{\cal U}}=\card{{\cal A}({\cal S})}=1$ for all ${\cal S}\subseteq{\cal N}$ with the balancedness condition $\displaystyle\sum_{{\cal S}\subsetneq{\cal N}}\mu({\cal S})v({\cal S})\leq v({\cal N})$. We can now state the following theorem regarding the existence of core decisions.

%\pagebreak
\begin{theorem}
\label{thm:corecond}
A robust cooperative game $({\cal N},{\cal A},{\cal V}({\cal U}))$ has a non-empty core if and only if it is balanced.
\end{theorem}
In order to prove Theorem~\ref{thm:corecond}, we need the following lemma.
\begin{lemma}
\label{lem:breakcond}
Given a robust cooperative game $({\cal N},{\cal A},{\cal V}({\cal U}))$ and an imputation $(a,\bz)$, a coalition ${\cal S}\subsetneq{\cal N}$ has the incentive to break away if and only if
\be
\label{eq:breakcond}
\sum_{i\in{\cal S}}z_i< v_{\max}(a,{\cal S}).
\ee
\end{lemma}

\begin{pf}
Given an imputation $(a,\bz)$, a coalition ${\cal S}$ has the incentive to break away if there exists an efficient decision $(\hat{a},\hat{\bz})$ %such that $\mb{{\cal X}}_{{\cal S}}({\cal U})\not\succsim_n\hat{\mb{{\cal X}}}_{{\cal S}}({\cal U})$. It means there exists $u\in{\cal U}$ 
and a realization $u\in{\cal U}$ such that for all $i\in{\cal S}$,
$$
v_u(\hat{a},{\cal S})\cdot \hat{z}_i>v_u(a,{\cal N})\cdot z_i.
$$
Since $(\hat{a},\hat{\bz})$ is efficient, we have: $\displaystyle\sum_{i\in{\cal S}}\hat{z}_i=1$. Summing over all $i\in{\cal S}$ the inequality above, we then obtain the following statement:
$$
\exists\,u\in{\cal U}\,:\,v_u(a,{\cal N})\cdot\sum_{i\in{\cal S}} z_i<v_u(\hat{a},{\cal S})\cdot \sum_{i\in{\cal S}}\hat{z}_i=v_u(\hat{a},{\cal S}).
$$
Equivalently, under Assumption \ref{as:pos}(ii), we have:
$$
\exists\,u\in{\cal U}\,:\,\sum_{i\in{\cal S}}z_i<\frac{v_u(\hat{a},{\cal S})}{v_u(a,{\cal N})}\,\Leftrightarrow\,\sum_{i\in{\cal S}}z_i<\max_{u\in{\cal U}}\left\{\frac{v_u(\hat{a},{\cal S})}{v_u(a,{\cal N})}\right\}.
$$
We have:
$$
\ba{rcl}
\displaystyle\max_{u\in{\cal U}}\left\{\frac{v_u(\hat{a},{\cal S})}{v_u(a,{\cal N})}\right\}&\leq&\displaystyle\max_{\alpha_{\cal S}\in{\cal A}({\cal S})}\max_{u\in{\cal U}}\left\{\frac{v_u(\alpha_{\cal S},{\cal S})}{v_u(a,{\cal N})}\right\}\\
&=&\displaystyle\max_{u\in{\cal U}}\max_{\alpha_{\cal S}\in{\cal A}({\cal S})}\left\{\frac{v_u(\alpha_{\cal S},{\cal S})}{v_u(a,{\cal N})}\right\}\\
&=&\displaystyle\max_{u\in{\cal U}}\left\{\frac{\displaystyle\max_{\alpha_{\cal S}\in{\cal A}({\cal S})}v_u(\alpha_{\cal S},{\cal S})}{v_u(a,{\cal N})}\right\}.
\ea
$$
The second equality is due to the fact that $v_u(a,{\cal N})>0$ for all $u\in{\cal U}$ under Assumption \ref{as:pos}(ii). Thus, if a coalition ${\cal S}$ has the incentive to break away, then
$$
\sum_{i\in{\cal S}}z_i<\max_{u\in{\cal U}}\left\{\frac{\displaystyle\max_{\alpha_{{\cal S}}\in{\cal A}({\cal S})}v_u(\alpha_{{\cal S}},{\cal S})}{v_u(a,{\cal N})}\right\}=v_{\max}(a,{\cal S}).
$$

Now, suppose $\displaystyle\sum_{i\in{\cal S}}z_i<v_{\max}(a,{\cal S})$, we will show that coalition ${\cal S}$ has the incentive to break away. We would need to show the existence of an efficient decision $(\hat{a},\hat{\bz})$ and a realization $u\in{\cal U}$ such that for all $i\in{\cal S}$, %such that $\mb{{\cal X}}_{{\cal S}}({\cal U})\not\succsim_n\hat{\mb{{\cal X}}}_{{\cal S}}({\cal U})$, i.e., there exists $u\in{\cal U}$ such that for all $i\in{\cal S}$,
$$
v_u(\hat{a},{\cal S})\cdot \hat{z}_i>v_u(a,{\cal N})\cdot z_i.
$$
Let $(\hat{a},\hat{u})\in\displaystyle\arg\max_{\alpha_{{\cal S}}\in{\cal A}({\cal S})}\max_{u\in{\cal U}}\left\{\frac{v_u(\alpha_{\cal S},{\cal S})}{v_u(a,{\cal N})}\right\}$. Clearly, we have:
$$
\frac{v_{\hat{u}}(\hat{a},{\cal S})}{v_{\hat{u}}(a,{\cal N})}=\max_{\alpha_{{\cal S}}\in{\cal A}({\cal S})}\max_{u\in{\cal U}}\left\{\frac{v_u(\alpha_{{\cal S}},{\cal S})}{v_u(a,{\cal N})}\right\}=\max_{u\in{\cal U}}\left\{\frac{\displaystyle\max_{\alpha_{\cal S}\in{\cal A}({\cal S})}v_u(\alpha_{\cal S},{\cal S})}{v_u(a,{\cal N})}\right\}=v_{\max}(a,{\cal S}).
$$
From Lemma~\ref{prop:individual:rationality:cond}, we have $\displaystyle z_i \geq v_{\max}(a,\{i\})$ since $(a,\bz)$ is an imputation. We also have, by Assumption~\ref{as:pos}(i), $\displaystyle v_{\max}(a,\{i\}) \geq 0$. Thus, $\displaystyle v_{\max}(a,{\cal S}) > \displaystyle\sum_{i\in{\cal S}}z_i \geq 0$ and hence $v_{\hat{u}}(\hat{a},{\cal S})\neq 0$.

Let $\displaystyle\eps=\frac{1}{\card{{\cal S}}}\left(v_{\max}(a,{\cal S})-\sum_{i\in{\cal S}}z_i\right)>0$ and define $\displaystyle\hat{z}_i=\frac{z_i+\eps}{v_{\max}(a,{\cal S})}$ for all $i\in{\cal S}$. Clearly, $\displaystyle\sum_{i\in{\cal S}}\hat{z}_i=1$ given the definition of $\eps$. In addition, for all $i\in{\cal S}$, we have:
$$
v_{\hat{u}}(\hat{a},{\cal S})\cdot\hat{z}_i=v_{\hat{u}}(a,{\cal N})\cdot z_i+v_{\hat{u}}(a,{\cal N})\cdot \eps>v_{\hat{u}}(a,{\cal N})\cdot z_i.
$$
Thus, $(\hat{a},\hat{\bz})$ is an efficient decision and the inequality above shows that coalition $\cal S$ indeed has the incentive to break away.
\end{pf}

We are now ready to prove Theorem \ref{thm:corecond}.

\begin{pf}
Suppose there exists a core decision $(a,\bz)$. Clearly, $(a,\bz)$ is an imputation. Applying Lemma \ref{lem:breakcond}, we can show that $\bz$ is a feasible (and optimal) solution of the following linear program:
\be
\label{eq:primal}
\ba{rl}
\displaystyle Z_P=\min_{\bx} & \displaystyle\sum_{i\in{\cal N}}0\cdot x_i\\
\st & \displaystyle\sum_{i\in{\cal S}}x_i\geq v_{\max}(a,{\cal S}),\quad \forall\,{\cal S}\subsetneq{\cal N},\\
& \displaystyle\sum_{i\in{\cal N}}x_i=1.
\ea
\ee
The dual problem is written as follows:
\be
\label{eq:dual}
\ba{rll}
Z_D=\max & \displaystyle\sum_{{\cal S}\subsetneq{\cal N}}y_{{\cal S}}v_{\max}(a,{\cal S})-p\\
\st & \displaystyle\sum_{{\cal S}:i\in{\cal S}}y_{{\cal S}}-p=0, & \forall\,i\in{\cal N},\\
& y_{{\cal S}}\geq 0, & \forall\,{\cal S}\subsetneq{\cal N}.
\ea
\ee
Applying strong duality, we have: $Z_D=Z_P=0$, which means, for all feasible solution $\left(\left\{y_{{\cal S}}\right\}_{{\cal S}\subsetneq{\cal N}},p\right)$,
$$
\sum_{{\cal S}\subsetneq{\cal N}}y_{{\cal S}}v_{\max}(a,{\cal S})\leq p.
$$
Now consider any balanced map $\mu$, it is clear that $\left(\left\{\mu({\cal S})\right\}_{{\cal S}\subsetneq{\cal N}},1\right)$ is a feasible solution of the dual problem. Thus we have:
$$
\sum_{{\cal S}\subsetneq{\cal N}}\mu({\cal S})v_{\max}(a,{\cal S})\leq 1.
$$
It shows that the robust cooperative game is balanced.

Now, suppose the robust cooperative game is balanced. There exists an action $a\in{\cal A}({\cal N})$ such that for all balanced map $\mu$,
$$
\sum_{{\cal S}\subsetneq{\cal N}}\mu({\cal S})v_{\max}(a,{\cal S})\leq 1.
$$
We are going to show that $Z_D=0$. The dual problem is indeed feasible. $(\mb{0},0)$ is a feasible solution and it implies that $Z_D\geq 0$. Let us consider a feasible solution $\left(\left\{y_{{\cal S}}\right\}_{{\cal S}\subsetneq{\cal N}},p\right)$. Since $y_{{\cal S}}\geq 0$ for all ${\cal S}\subsetneq{\cal N}$, we have: $p\geq 0$. If $p=0$, it is easy to show that $y_{{\cal S}}=0$ for all ${\cal S}\subsetneq{\cal N}$. If $p>0$, let $\mu({\cal S})=y_{\cal S}/p$ for all ${\cal S}\subsetneq{\cal N}$. The map $\mu$ is indeed a balanced map, thus we have:
$$
\sum_{{\cal S}\subsetneq{\cal N}}\left(y_{{\cal S}}/p\right)v_{\max}(a,{\cal S})\leq 1\,\Leftrightarrow\,\sum_{{\cal S}\subsetneq{\cal N}}y_{{\cal S}}v_{\max}(a,{\cal S})-p\leq 0.
$$
Thus for all feasible solution $\left(\left\{y_{{\cal S}}\right\}_{{\cal S}\subsetneq{\cal N}},p\right)$, $\displaystyle \sum_{{\cal S}\subsetneq{\cal N}}y_{{\cal S}}v_{\max}(a,{\cal S})-p\leq 0$, which means $Z_D\leq 0$ and hence $Z_D=0$. According linear strong duality, the primal problem is feasible (and with the obvious optimal objective $Z_P=0$). Thus there exists at least a feasible solution $\bz$ of the primal problem. The decision $(a,\bz)$ is an imputation and indeed a core decision given the conditions obtained from the constraints of the primal problem. Thus the robust cooperative game has a non-empty core.
\end{pf}

Theorem \refs{thm:corecond} establishes the relationship between the game balancedness and the existence of core solutions. Computationally, it implies that we can attempt to solve the following optimization problem to check the existence of core solutions:
\be
\label{eq:rleastcore}
\ba{rl}
\displaystyle s({\cal N},{\cal A},{\cal V}({\cal U}))=\min_{\bx,\eps,a} & \eps\\
\st & \displaystyle\sum_{i\in{\cal S}}x_i\geq v_{\max}(a,{\cal S})-\eps,\quad \forall\,{\cal S}\subsetneq{\cal N},\\
& \displaystyle\sum_{i\in{\cal N}}x_i=1,\\
& a\in{\cal A}({\cal N}).
\ea
\ee
If $s({\cal N},{\cal A},{\cal V}({\cal U}))\leq 0$, then the core of the robust cooperative game $({\cal N},{\cal A},{\cal V}({\cal U}))$ is non-empty. Note that this optimization problem is no longer a linear program like Problem \refs{eq:leastcore} or \refs{eq:stability} given the fact that the action of the grand coalition is now a decision variable. We will investigate further the computational aspect of finding a core solution for a specific robust cooperative game, the \emph{robust newsvendor game}, which is going to be discussed next. Even though the framework of robust cooperative games developed in this section is indeed suitable for the newsvendor games with distributional ambiguity which we are interested in, we would like to emphasize that it is plausible to develop other frameworks of cooperative games with uncertain characteristic functions using different payoff distribution schemes and different preference relations to suit other applications better.
%We will focus on an in-deep study of the newsvendor games with ambiguity in demand distribution in the next section.

%\section{Stochastic Newsvendor Games}
%\label{sec:newsvendor}

\section{Robust Newsvendor Games} \label{stoc.newsvendor}
In this section, we consider newsvendor games with ambiguity in demand distributions in the framework of robust cooperative games, which we call \emph{robust newsvendor games}.
\begin{comment}Individual retailers usually collect historical demands independently before they join any coalition and the current assumption of a known joint distribution can be considered quite strong. In order to make the problem is more realistic, we only assume that some (multivariate) marginal distributions are known. For example, it is more reasonable to assume the knowledge of the joint demand distribution of a subset of retailers that are located close to each other and hence likely serve customers from a same area.
\end{comment}
As discussed in the previous section, the uncertainty comes from the fact that the joint demand distribution is unknown. We assume that only some (multivariate) marginal distributions of the joint demand are known. More concretely, consider a partition of ${\cal N}$ with $R$ subsets ${\cal N}_1,\ldots,{\cal N}_R$ such that
$$
{\cal N}=\bigcup_{r=1}^R{\cal N}_r\quad\mbox{and}\quad{\cal N}_r\cap{\cal N}_s=\emptyset\quad\mbox{for all }r\neq s.
$$
Given a vector $\mb{d}\in\mathbb{R}^n$, let $\mb{d}_r\in\mathbb{R}^{N_r}$ denote the sub-vector formed with the elements in the $r$th subset ${\cal N}_r$ where $N_r=\card{{\cal N}_r}$ is the size of the subset. We assume that  probability measures $P_r$ of random vectors $\td{\mb{d}}_r$ are known for all $r=1,\ldots,R$. Let $\mathcal{P}(P_1,\ldots,P_R)$ denote the set of joint probability measures of the random vector $\td{\mb{d}}$ consistent with the prescribed probability measures of the random vectors $\td{\mb{d}}_r$ for all $r=1,\ldots,R$, which acts as the uncertainty set ${\cal U}$ in the general framework of robust cooperative games. Note that $\mathcal{P}(P_1,\ldots,P_R)$ is always non-empty since the independent measure among the sub-vectors is a feasible distribution. The set of joint distributions with fixed marginal distributions $\mathcal{P}(P_1,\ldots,P_R)$ is referred to as the Fr\'echet class of distributions (see R\"uschendorf \cite{ruschen91b}). It has been used to evaluate bounds on the cumulative distribution function of a sum of random variables with an application in risk management (Embrechts and Puccetti \cite{ep06b}). Doan and Natarajan \cite{doan12} developed a robust optimization model using $\mathcal{P}(P_1,\ldots,P_R)$ with an application in project management. We now investigate this Fr\'echet class of distributions in the context of robust newsvendor games.

Given a subset ${\cal S}\subseteq{\cal N}$, we define ${\cal S}_r={\cal S}\cap{\cal N}_r$ for all $r=1,\dots,R$. Clearly, if all retailers $i$, $i\in{\cal S}$, join together, we know the non-overlapping marginal distributions of the joint demand vector with respect to the partition $({\cal S}_1,\ldots,{\cal S}_R)$ of ${\cal S}$. If ${\cal S}\subseteq{\cal N}_r$ for some $r$, the joint distribution of $\td{d}_i$, $i\in{\cal S}$, is completely known. In this case, the action the coalition should take, i.e., the decision on the ordering quantity, is well-defined as in the deterministic setting. The action set of coalition ${\cal S}$ can be simply defined as ${\cal Y}({\cal S})=\{y^*({\cal S})\}$, where $y^*({\cal S})$ is the $(p-c)/p$-quantile of the known distribution of $\tilde{d}({\cal S})$ as defined in \refs{eq:optquantbar}. Note that ${\cal Y}({\cal S})$ acts as the action set ${\cal A}({\cal S})$ in the general framework of robust cooperative games. In general, the distribution $P({\cal S})$ is unknown and coalition ${\cal S}$ can choose its action regarding the ordering quantity from a general action set ${\cal Y}({\cal S})$. %If the random demand $\tilde{d}({\cal S})$ is bounded, i.e., $\td{d}({\cal S})\in[d_{\min}({\cal S});d_{\max}({\cal S})]$, the action set could be set as ${\cal Y}({\cal S})=[d_{\min}({\cal S});d_{\max}({\cal S})]$ given that it is not profitable to order less than the minimum quantity $d_{\min}({\cal S})$ or more than the maximum quantity $d_{\max}({\cal S})$. (Indeed, for $0\leq y\leq d_{\min}({\cal S})$, $v(y,{\cal S})=(p-c)y$, which is a strictly increasing function in $y$ given that $p>c$. Similarly, for $y\geq d_{\max}({\cal S})$, $v(y,{\cal S})=-cy + p\,\mathbb{E}_{P({\cal S})}\left[\td{d}({\cal S})\right]$, which is a strictly decreasing function in $y$ given $c>0$.) Having said that, %
To keep it simple, we shall let ${\cal Y}({\cal S})=\mathbb{R}_+$ given the fact that the ordering quantities are non-negative for all ${\cal S}\subseteq{\cal N}$. %, except for all subsets of ${\cal N}_r$, $r=1,\ldots,R,$ and the grand coalition ${\cal N}$.
 As mentioned earlier, if ${\cal S}\subseteq{\cal N}_r$ for some $r$, we can restrict ${\cal Y}({\cal S})=\{y^*({\cal S})\}$, where $y^*({\cal S})$ is the $(p-c)/p$-quantile of the known distribution of $\tilde{d}({\cal S})$ as defined in \refs{eq:optquantbar}.

Given an ordering decision $y\in{\cal Y}({\cal S})$, the \emph{uncertain} payoff, i.e., total expected profit, of coalition $\cal S$ is $v_P(y,{\cal S})=\mathbb{E}_{P({\cal S})}\left[p\min\{\td{d}({\cal S}),y\}-cy\right]$ for $P\in{\cal P}(P_1,\ldots,P_R)$ as in \refs{eq:quantbar}, where $P({\cal S})$ is the corresponding marginal joint distribution of $\td{d}_i$, $i\in{\cal S}$, derived from $P$. For each individual retailer $i$, ${\cal Y}(\{i\})=y_i^*$ and $v_P(y_i^*,\{i\})\geq v_P(0,\{i\})=0$ for all $P\in{\cal P}(P_1,\ldots,P_R)$, which implies Assumption~\ref{as:pos}(i) is automatically satisfied.

Finally, for the grand coalition, in order to satisfy Assumption~\ref{as:pos}(ii), we let
\be
\label{eq:yn}
{\cal Y}({\cal N})=\left\{y\in\mathbb{R}_+\,:\,v_P(y,{\cal N})>0,\,\forall\,P\in{\cal P}(P_1,\ldots,P_R) \right\}.
\ee
We shall provide a simple condition with which the action set of the grand coalition is non-empty. Let $d_{\min}({\cal S})$ be the minimum value that the random demand $\td{d}({\cal S})$ can achieve, the following lemma sets out a sufficient condition for ${\cal Y}({\cal N})$ to be empty.

\begin{lemma}
\label{lem:nonempty}
If $\displaystyle\max_{r=1,\ldots,R}d_{\min}({\cal N}_r)>0$, then ${\cal Y}({\cal N})\neq \emptyset$.
\end{lemma}

\begin{pf}
We have: $\displaystyle\td{d}({\cal N})=\sum_{r=1}^R\td{d}({\cal N}_r)$. Thus, if $\displaystyle\max_{r=1,\ldots,R}d_{\min}({\cal N}_r)>0$, then:
$$
d_{\min}({\cal N})\geq\sum_{r=1}^Rd_{\min}({\cal N}_r)> 0.
$$
For $y\in[0,d_{\min}(N)]$, $v_P(y,{\cal N})=(p-c)y$, which is strictly increasing for any $P\in{\cal P}(P_1,\ldots,P_R)$ given the fact that $p>c$. In addition, for $y\geq d_{\max}({\cal N})$, $v_P(y,{\cal N})=-cy+\displaystyle p\sum_{r=1}^R\mathbb{E}_{P_r}\left[\td{d}({\cal N}_r)\right]$, which is strictly decreasing for any $P\in{\cal P}(P_1,\ldots,P_R)$ given the fact that $c>0$.

Now consider the function $\displaystyle\bar{v}(y,{\cal N})=\min_{P\in{\cal P}(P_1,\ldots,P_R)}v_P(y,{\cal N})$. It is clear that $\bar{v}(\cdot,{\cal N})$ is again strictly increasing in $[0,d_{\min}({\cal N})]$ and strictly decreasing in $[d_{\max}({\cal N}),+\infty)$. Thus we have:
$$
\arg\max_{y\geq 0}\bar{v}(y,{\cal N})\in[d_{\min}({\cal N}),d_{\max}({\cal N})],
$$
and
$$
\max_{y\geq 0}\bar{v}(y,{\cal N})\geq \bar{v}(d_{\min}({\cal N}),{\cal N})=(p-c)d_{\min}({\cal N})>0.
$$
Thus there exists $y\geq 0$ such that $\bar{v}(y,{\cal N})>0$, or equivalently, $v_P(y,{\cal N})>0$ for all $P\in{\cal P}(P_1,\ldots,P_R)$. It shows that ${\cal Y}({\cal N})\neq\emptyset$.
\end{pf}

The condition in Lemma \ref{lem:nonempty} simply requires that one of the (marginal) distributions $P_r$, $r=1,\ldots,R$ has the support set of solely non-zero demand vectors, which can be considered as a reasonable assumption in reality. %The proof of Lemma \ref{lem:nonempty} also indicates that ${\cal Y}({\cal N})$ is an interval if it is not empty.
 In the rest of the paper, we shall make that assumption to ensure ${\cal Y}({\cal N})\neq\emptyset$, or equivalently, that Assumption~\ref{as:pos}(ii) is satisfied.
% \begin{assumption}
% \label{as:mpos}
% There exists $r=1,\ldots,R$ such that $d_{\min}({\cal N}_r)>0$.
% \end{assumption}
We are now ready to consider the robust newsvendor game $({\cal N},{\cal Y},{\cal V}({\cal P}(P_1,\ldots,P_r)))$ and investigate the existence of its imputations and core solutions.
\subsection{Existence of Imputations and Core Solutions}
\label{ssec:existence}
The uncertain characteristic function $v_P(y,{\cal S})$ can be written as
\be
\label{eq:char}
v_{P}(y,{\cal S})=(p-c)y-p\,\mathbb{E}_{P({\cal S})}\left[\left(y-\td{d}({\cal S})\right)^+\right],
\ee
for all ${\cal S}\subseteq{\cal N}$, $y\in{\cal Y}({\cal S})$, and $P\in{\cal P}(P_1,\ldots,P_R)$. The deterministic newsvendor games always have imputations since the corresponding characteristic function is super-additive. The following theorem claims the existence of imputations of the robust newsvendor game $({\cal N},{\cal Y},{\cal V}({\cal P}(P_1,\ldots,P_r)))$.
\begin{theorem}
\label{thm:rimpu}
The robust newsvendor game $({\cal N},{\cal Y},{\cal V}({\cal P}(P_1,\ldots,P_r)))$ always has an imputation.
\end{theorem}

In order to prove Theorem \ref{thm:rimpu}, we first study a particular action that each coalition can take, the worst-case optimal ordering quantity $y^*_{wc}({\cal S})$:
\be
\label{eq:optquant2}
y^*_{wc}({\cal S})\in\arg\max_{y\geq 0}\left\{(p-c)y-p\max_{P\in{\cal P}(P_1,\ldots,P_R)}\mathbb{E}_{P}\left[\left(y-\td{{d}}({\cal S})\right)^+\right]\right\}.
\ee
Basically, this ordering quantity is optimal under the worst-case scenario with respect to the joint demand distribution $P$. Let us also define $v_{wc}({\cal S})$ as the maximum worst-case expected profit for coalition ${\cal S}$, i.e.
\be
v_{wc}({\cal S}) = \max_{y\geq 0}\left\{(p-c)y-p\max_{P\in{\cal P}(P_1,\ldots,P_R)}\mathbb{E}_{P}\left[\left(y-\td{{d}}({\cal S})\right)^+\right]\right\}.
\ee

When ${\cal S}\subseteq{\cal N}_r$ for some $r$, clearly, $y^*_{wc}({\cal S})=y^*({\cal S})$, the $(p-c)/p$-quantile of the known distribution of $\td{d}({\cal S})$, and $v_{wc}({\cal S})=\bar{v}({\cal S})$ as defined in \refs{eq:optval}. The following lemma shows how to calculate the worst-case optimal ordering quantities $y^*_{wc}({\cal S})$ and the worst-case expected profit $v_{wc}({\cal S})$ for an arbitrary ${\cal S}\subseteq{\cal N}$.
\begin{lemma}
\label{lem:optquant}
For an arbitrary ${\cal S}\subseteq{\cal N}$, the worst-case optimal ordering quantity $y^*_{wc}({\cal S})$ defined in \refs{eq:optquant2} can be calculated as follows:
\be
\label{eq:optquantform1}
y^*_{wc}({\cal S})=\sum_{r=1}^R{y}^*({\cal S}_r),
\ee
where ${\cal S}_r={\cal S}\cap{\cal N}_r$ for all $r=1,\ldots,R$, ${y}^*({\cal S}_r)$ is the $(p-c)/p$-quantile of the known distribution of $\td{d}({\cal S}_r)$, and ${y}^*(\emptyset)=0$. In addition, $v_{wc}({\cal S}) = \displaystyle \sum_{r=1}^R v_{wc}({\cal S}_r) = \displaystyle \sum_{r=1}^R \bar{v}({\cal S}_r)$.
\end{lemma}

\begin{pf}
Consider the optimization problem in \refs{eq:optquant2}. For $y\leq d_{\min}({\cal S})=\min\{\td{d}({\cal S})\}$, we can write $(p-c)y-p\,\mathbb{E}_{P}\left[\left(y-\td{{d}}({\cal S})\right)^+\right]=(p-c)y$ for any distribution $P$. Since $p-c>0$, it is an increasing function in $y$ in $(-\infty;d_{\min}({\cal S})]$. Since $d_{\min}({\cal S})\geq 0$, we can then remove the non-negative constraint $y\geq 0$ from \refs{eq:optquant2} when calculating $y^*_{wc}({\cal S})$. Now, consider the inner optimization problem of \refs{eq:optquant2}. This is an instance of the distributionally robust optimization problem studied in Doan and Natarajan \cite{doan12}. Without loss of generality, we can assume that ${\cal S}_r\neq\emptyset$ for all $r=1,\ldots,R$ knowing that ${y}^*(\emptyset)=0$. Applying Proposition 1(ii) from \cite{doan12}, we obtain the following reformulation:
$$
\ba{rl}
\displaystyle\max_{P\in{\cal P}(P_1,\ldots,P_R)}\mathbb{E}_{P}\left[\left(y-\td{{d}}({\cal S})\right)^+\right]=\min_{\mbs{x}} &\displaystyle\sum_{r=1}^R\mathbb{E}_{P_r}\left[\left(x_r-\td{d}({\cal S}_r)\right)^+\right]\\
\st &\displaystyle\sum_{r=1}^Rx_r=y.
\ea
$$
Thus, in order to find $y^*_{wc}({\cal S})$, we can solve the following optimization problem
$$
\ba{rl}
\displaystyle\max_{y,\mbs{x}}&\displaystyle (p-c)y-p\sum_{r=1}^R\mathbb{E}_{P_r}\left[\left(x_r-\td{d}({\cal S}_r)\right)^+\right]\\
\st & \displaystyle\sum_{r=1}^Rx_r=y.\\
\ea
$$
The optimal ordering quantity $y^*_{wc}({\cal S})$ can then be calculated as $\displaystyle y^*_{wc}({\cal S})=\sum_{r=1}^Rx_r^*$, where $\bx^*$ is the optimal solution of the following separable optimization problem:
$$
\max_{\mbs{x}}\sum_{r=1}^R\left((p-c)x_r-p\sum_{r=1}^R\mathbb{E}_{P_r}\left[\left(x_r-\td{d}({\cal S}_r)\right)^+\right]\right),$$
or equivalently,
$$
\sum_{r=1}^R\max_{x_r}\left((p-c)x_r-p\mathbb{E}_{P_r}\left[\left(x_r-\td{d}({\cal S}_r)\right)^+\right]\right).
$$
For each sub-problem, $x_r^*$ is the $(p-c)/p$-quantile of the distribution of $\td{d}({\cal S}_r)$, which means $x^*_r={y}^*({\cal S}_r)$ for all $r=1,\ldots, R$, according to \refs{eq:optquantbar}. Thus we have $\displaystyle y^*_{wc}({\cal S})=\sum_{r=1}^R{y}^*({\cal S}_r)$. This also leads to $v_{wc}({\cal S}) = \displaystyle \sum_{r=1}^R v_{wc}({\cal S}_r) = \displaystyle \sum_{r=1}^R \bar{v}({\cal S}_r)$. Here, the joint demand distribution for players in ${\cal S}_r$ is known and hence the worst-case expected profit ${v}_{wc}(S_r)$ is exactly the same with the deterministic expected profit $\bar{v}({\cal S}_r)$ shown in Formulation~\refs{eq:optval}.
\end{pf}

%Lemma \ref{lem:optquant} shows that the robust optimal ordering quantity $y^*_{wc}({\cal S})$ can be computed using the optimal ordering quantities $y^*({\cal S}_r)$, which is the $(p-c)/p$-quantile of the known distribution of $\td{d}({\cal S}_r)$, where ${\cal S}_r={\cal S}\cup{\cal N}_r$ for $r=1,\ldots, R$. All of these quantities $y^*({\cal S})$, where ${\cal S}\subseteq{\cal N}_r$, $r=1,\ldots,R$, can be easily computed given the fact that the marginal distributions $P_1,\ldots,P_R$, are completely known.
 We are now ready to use this lemma to prove Theorem \ref{thm:rimpu}.

\begin{pf}
According to Theorem \ref{prop:impucond}, the robust newsvendor game $({\cal N},{\cal Y},{\cal V}({\cal P}(P_1,\ldots,P_r)))$ has imputations if and only if there exists $y\in{\cal Y}({\cal N})$ such that
$$
\sum_{i\in{\cal N}}v_{\max}(y,\{i\})\leq 1,
$$
where $\displaystyle v_{\max}(y,\{i\})=\max_{P\in{\cal P}(P_1,\ldots,P_R)}\left\{\frac{\displaystyle\max_{y_i\in{\cal Y}(\{i\})}v_P(y_i,\{i\})}{v_P(y,{\cal N})}\right\}$. For each individual retailer $i$, the demand distribution is known; therefore ${\cal Y}(\{i\})=\{y_i^*\}$, where $y_i^*$ is the $(p-c)/p$-quantile of the distribution function of $\td{d}_i$. Thus
$$
\max_{y_i\in{\cal Y}(\{i\})}v_P(y_i,\{i\})=v_i(y^*)=\bar{v}_i\geq 0,\quad\forall\,P\in{\cal P}(P_1,\ldots,P_R).
$$
We then have $\displaystyle v_{\max}(y,\{i\})=\bar{v}_i\cdot\left(\displaystyle\min_{P\in{\cal P}(P_1,\ldots,P_R)}v_P(y,{\cal N})\right)^{-1}$ and the existence condition of imputation can be written as follows:
$$
\exists\,y\in{\cal Y}({\cal N}):\,\sum_{i\in{\cal N}}\bar{v}_i\leq \displaystyle\min_{P\in{\cal P}(P_1,\ldots,P_R)}v_P(y,{\cal N}),
$$
given the fact that $v_P(y,{\cal N})>0$ for all $y\in{\cal Y}({\cal N})$ and $P\in{\cal P}(P_1,\ldots,P_R)$. Equivalently, the existence condition is
$$
\sum_{i\in{\cal N}}\bar{v}_i\leq \displaystyle\max_{y\in{\cal Y}({\cal N})}\min_{P\in{\cal P}(P_1,\ldots,P_R)}v_P(y,{\cal N})=\max_{y\geq 0}\min_{P\in{\cal P}(P_1,\ldots,P_R)}v_P(y,{\cal N}).
$$
The optimization problem $\displaystyle\max_{y\geq 0}\min_{P\in{\cal P}(P_1,\ldots,P_R)}v_P(y,{\cal N})$ is the same as Problem \refs{eq:optquant2}. Let us consider the worst-case optimal ordering quantity $y^*_{wc}({\cal N})$, the existence condition can then be written as follows:
$$
\sum_{i\in{\cal N}}\bar{v}_i\leq v_P(y^*_{wc}({\cal N}),{\cal N}),\quad\forall\,P\in{\cal P}(P_1,\ldots,P_R).
$$
Applying Lemma \ref{lem:optquant} for ${\cal S}={\cal N}$, we have $\displaystyle y^*_{wc}({\cal N})=\sum_{i=1}^R{y}^*({\cal N}_r)$. Using the fact that $(x+y)^+\leq x^++y^+$, we have, for all $P\in{\cal P}(P_1,\ldots,P_R)$,
$$
\ba{rl}
v_P(y^*_{wc}({\cal N}),{\cal N}) &=\,(p-c)y^*_{wc}({\cal N})-p\,\mathbb{E}_P\left[\left(y^*_{wc}({\cal N})-\td{d}({\cal N})\right)^+\right]\\
&\geq\,\displaystyle (p-c)\sum_{i=1}^R{y}^*({\cal N}_r)-p\sum_{r=1}^R\mathbb{E}_{P_r}\left[\left({y}^*({\cal N}_r)-\td{d}({\cal N}_r)\right)^+\right]\\
&=\,\displaystyle\sum_{r=1}^R\left\{(p-c){y}^*({\cal N}_r)-p\,\mathbb{E}_{P_r}\left[\left({y}^*({\cal N}_r)-\td{d}({\cal N}_r)\right)^+\right]\right\}\\
&=\,\displaystyle\sum_{r=1}^R\bar{v}({\cal N}_r),
\ea
$$
where $\bar{v}({\cal N}_r)$ is the optimal total expected profit of coalition ${\cal N}_r$ and is computed using \refs{eq:optval} for all $r=1,\ldots,R$ since $P_1,\ldots,P_R$ are completely known.
\begin{comment}
We have $\bar{v}({\cal N}_r)\geq 0$ for all $r=1,\ldots,R$ and under Assumption \ref{as:mpos}, there exists $r$ such that $\bar{v}({\cal N}_r)>0$. Thus, we have
$$
v_P({\cal N})\geq\sum_{r=1}^R\bar{v}({\cal N}_r)>0.
$$
\end{comment}

Now consider the deterministic newsvendor game for coalition ${\cal N}_r$ with the complete knowledge of the joint distribution $P_r$. There exists at least one imputation for this cooperative game; thus, we have
$$
\sum_{i\in{\cal N}_r}\bar{v}_i\leq\bar{v}({\cal N}_r),\quad\forall\,r=1,\ldots,R.
$$
Using the fact that ${\cal N}=\bigcup_{r=1}^R{\cal N}_r$ and ${\cal N}_r\cap{\cal N}_s = \emptyset$ for all $r\neq s$, we have
$$
v_P({\cal N})\geq\sum_{r=1}^R\bar{v}({\cal N}_r)\geq\sum_{r=1}^R\sum_{i\in{\cal N}_r}\bar{v}_i=\sum_{i\in{\cal N}}\bar{v}_i,\quad\forall\,P\in{\cal P}(P_1,\ldots,P_R).
$$
Thus, the existence condition is satisfied and the robust newsvendor game $({\cal N},{\cal Y},{\cal V}({\cal P}(P_1,\ldots,P_r)))$ always has an imputation.
\end{pf}

We focus on properties of core solutions of the robust newsvendor game  $({\cal N},{\cal Y},{\cal V}({\cal P}(P_1,\ldots,P_r)))$ next.
\begin{theorem}
\label{thm:rcore2}
If $(y,\mb{z})$ is a core solution of the robust newsvendor game $({\cal N},{\cal Y},{\cal V}({\cal P}(P_1,\ldots,P_r)))$, then the following statements hold:
\begin{itemize}
\item[(a)] The ordering quantity $y$ is the worst-case optimal ordering quantity, i.e., $y = y^*_{wc}({\cal N})$.
\item[(b)] $v_{wc}({\cal N})\cdot\mb{z}({\cal N}_r)$ is a core solution of the (deterministic) newsvendor game $({\cal N}_r,\bar{v})$. %, i.e.
% \begin{itemize}
% \item[(b1)] $\displaystyle \left(\sum_{j=1}^R \bar{v}({\cal N}_j)\right) \sum_{i \in N_r} z_i = \bar{v}({\cal N}_r)$ and
% \item[(b2)] For all coalition ${\cal S}_r \subset {\cal N}_r$, we have $\displaystyle \left(\sum_{j=1}^R \bar{v}({\cal N}_j)\right) \sum_{i \in S_r} z_i \geq \bar{v}({\cal S}_r)$.
% \end{itemize}
\end{itemize}
\end{theorem}

\begin{pf}
Let us consider a core solution $(y,\mb{z})$ of the robust newsvendor game {\small $({\cal N},{\cal Y},{\cal V}({\cal P}(P_1,\ldots,P_r)))$}. For each coalition ${\cal N}_r$, $r =1,\ldots,R$, we have

%$\displaystyle v_{\max}(a,{\cal S})=\max_{u\in{\cal U}}\left\{\frac{\displaystyle\max_{\alpha_{{\cal S}}\in{\cal A}({\cal S})}v_u(\alpha_{{\cal S}},{\cal S})}{v_u(a,{\cal N})}\right\}$.
\begin{eqnarray}
v_{\max}(y,{\cal N}_r)&=&\max_{P\in{\cal P}(P_1,\ldots,P_R)} \left\{\frac{\displaystyle\max_{y_r \in{\cal Y}({\cal N}_r)}v_P(y_r,{\cal N}_r)}{v_P(y,{\cal N})}\right\} \nonumber\\
&=& \frac{\displaystyle\max_{y_r \in{\cal Y}({\cal N}_r)}v_{P_r}(y_r,{\cal N}_r)}{\displaystyle \min_{P\in{\cal P}(P_1,\ldots,P_R)} v_P(y,{\cal N})} \nonumber\\
&=& \frac{\displaystyle v_{wc}({\cal N}_r)}{\displaystyle \min_{P\in{\cal P}(P_1,\ldots,P_R)} v_P(y,{\cal N})} = \frac{\displaystyle \bar{v}({\cal N}_r)}{\displaystyle \min_{P\in{\cal P}(P_1,\ldots,P_R)} v_P(y,{\cal N})}. \label{eq:thm:rcore1}
\end{eqnarray}
This is due to the fact that $P_r$ is known with certainty and $v_{wc}({\cal N}_r)\geq 0$. According to Lemma~\ref{lem:breakcond} and Equality \refs{eq:thm:rcore1}, for $(y,\mb{z})$ to be a core solution, we need to have
\begin{eqnarray}
\displaystyle \sum_{i\in{\cal N}_r}z_i \geq v_{\max}(y,{\cal N}_r) = \frac{\displaystyle \bar{v}({\cal N}_r)}{\displaystyle \min_{P\in{\cal P}(P_1,\ldots,P_R)} v_P(y,{\cal N})}. \label{eq:thm:rcore2}
\end{eqnarray}
Summing this over $r=1,\ldots,R,$ we obtain
\begin{eqnarray}
\displaystyle 1 &=& \sum_{i\in {\cal N}} z_i = \sum_{r=1}^R~\sum_{i\in{\cal N}_r}z_i \geq \frac{\displaystyle \sum_{r=1}^R \bar{v}({\cal N}_r)}{\displaystyle \min_{P\in{\cal P}(P_1,\ldots,P_R)} v_P(y,{\cal N})}.\label{eq:thm:rcore3}
\end{eqnarray}
Given Assumption~\ref{as:pos}(ii), $\displaystyle \min_{P\in{\cal P}(P_1,\ldots,P_R)} v_P(a,{\cal N})>0$. Hence, we can rewrite Inequality~\refs{eq:thm:rcore3} as
\begin{eqnarray*}
\displaystyle \sum_{r=1}^R \bar{v}({\cal N}_r) \leq \min_{P\in{\cal P}(P_1,\ldots,P_R)} v_P(y,{\cal N}).
\end{eqnarray*}
This leads to
\begin{eqnarray}
\displaystyle \sum_{r=1}^R \bar{v}({\cal N}_r) &\leq& \max_{y \geq 0}~\min_{P\in{\cal P}(P_1,\ldots,P_R)} v_P(y,{\cal N}) \label{eq:thm:rcore4} \\
&\equiv& v_{wc}({\cal N}) \nonumber\\
&=& \sum_{r=1}^R v_{wc}({\cal N}_r) = \sum_{r=1}^R \bar{v}({\cal N}_r), \nonumber
\end{eqnarray}
where the last equality comes from Lemma~\ref{lem:optquant}. Thus, all the inequalities in the chain need to be tight. It implies that the ordering quantity $y$ is the worst-case optimal ordering quantity, $y=y_{wc}^*({\cal N})$ for \refs{eq:thm:rcore4} to be tight. We also obtain $\displaystyle \sum_{i \in N_r} z_i = \displaystyle \frac{\bar{v}({\cal N}_r)}{v_{wc}({\cal N})}$ for \refs{eq:thm:rcore2} to be tight. In addition, for all coalition ${\cal S}_r \subset {\cal N}_r$, by using the similar argument in deriving Equality~\refs{eq:thm:rcore1}, we have
\begin{eqnarray*}
\displaystyle \sum_{i \in S_r} z_i \geq v_{\max}(y,{\cal S}_r) = \displaystyle \frac{\bar{v}(S_r)}{v_{wc}({\cal N})},
\end{eqnarray*}
which means $v_{wc}({\cal N})\cdot\bz({\cal N}_r)$ is a core solution of the deterministic newsvendor game $({\cal N}_r,\bar{v})$ for $r=1,\ldots,R$.
\end{pf}

Theorem~\ref{thm:rcore2} shows that in order to check whether a particular decision $(y,\bz)$ is a core solution of the robust newsvendor game $({\cal N},{\cal Y},{\cal V}({\cal P}(P_1,\ldots,P_r)))$, we only need to consider $y=y^*_{wc}({\cal N})$. The optimization problem \refs{eq:rleastcore} for checking the existence of core solutions can be reduced to the following linear program:
\be
\label{eq:rnleastcore}
\ba{rl}
\displaystyle \sigma(y^*_{wc}({\cal N}))=\min_{\bx,\eps} & \eps\\
\st & \displaystyle\sum_{i\in{\cal S}}x_i\geq v_{\max}(y^*_{wc}({\cal N}),{\cal S})-\eps,\quad{\cal S}\subsetneq{\cal N},\\
& \displaystyle\sum_{i\in{\cal N}}x_i=1.
\ea
\ee

%In addition, the allocation $\mb{z}({\cal N}_r)$ for each subgroup of players corresponds to a stable payoff sharing among the deterministic subgame on that subset of players.

Unlike the deterministic newsvendor games, the robust newsvendor games do not always have core solution for $N\geq 3$. The following example shows a simple robust newsvendor game with $N=3$ whose core is empty.

\begin{example}
\label{ex:nocore}
Let us consider the partition of ${\cal N}=\{1,2,3\}$ with $R=2$, ${\cal N}_1=\{1,2\}$ and ${\cal N}_3=\{3\}$. The probability distribution $P_1$ of $(\td{d}_1,\td{d}_2)$ is characterized by the uniform marginal distribution of $\td{d}_1$, $\td{d}_1 \sim U(0,D)$, for some $D>0$, and the relationship $\td{d}_2 = D-\td{d}_1$. The probability distribution $P_2$ of $\td{d}_3$ is another uniform distribution, $\td{d}_3 \sim U(0,D)$. Clearly, ${\cal Y}({\cal N})\neq\emptyset$ given the fact that $d_{\min}({\cal N}_1)=D>0$.

According Theorem \ref{thm:rcore2}, if the robust newsvendor game $({\cal N},{\cal Y},{\cal V}({\cal P}(P_1,P_2)))$ has a core solution $(y,\bz)$, then $y=y^*_{wc}({\cal N})$. Given $P_1$ and $P_2$, we are able to compute and bound some values of $v_{\max}(y,{\cal S})$ where ${\cal S}\subset{\cal N}$, as follows:
$$
v_{\max}(y,\{1,2\})=\frac{2p}{3p-c}\leq v_{\max}(y,\{2,3\})=v_{\max}(y,\{1,3\}).
$$
For clarity of the exposition, we leave the detailed computation of these values in the Appendix. Now, since $(y,\bz)$ is a core solution, $\displaystyle\sum_{i\in{\cal S}}z_i\geq v_{\max}(y,{\cal S})$ for ${\cal S}\subset{\cal N}$ and $z_1+z_2+z_3=1$. Thus
\begin{eqnarray*}
2 &=& (z_1+z_2)+(z_2+z_3)+(z_3+z_1)\\
&\geq& v_{max}(y,\{1,2\}) + v_{max}(y,\{2,3\})+v_{max}(y,\{1,3\})\\
&\geq& 6p/(3p-c) > 2,
\end{eqnarray*}
which is a contradiction or this robust newsvendor game does not have a core solution.
\end{example}

We now focus on how to solve Problem \refs{eq:rnleastcore} to check the existence of core solutions of robust newsvendor games. If the core is empty, i.e., $\sigma(y^*_{wc}({\cal N}))>0$, \emph{least core solutions} for the robust newsvendor game can be found by solving the general problem \refs{eq:rleastcore}, which will be discussed in the next section.
\subsection{Core and Least Core Computation}\label{robust_core_computation}
\begin{comment}
In order to solve the problem \refs{eq:rnleastcore}, we need to be able to compute $v_{\max}(y^*_{wc}({\cal N}),{\cal S})$ for each coalition $\cal S\subsetneq{\cal N}$. If ${\cal S}\subseteq {\cal N}_r$ for some $r$, we have shown that $\displaystyle v_{\max}(y^*_{wc}({\cal N}),{\cal S})=\frac{\bar{v}(S)}{v_{wc}({\cal N})}$.
\end{comment}
Both Problems \refs{eq:rleastcore} and \refs{eq:rnleastcore} involve the function $v_{\max}$. We first show how to compute $v_{\max}(y,{\cal S})$ for an arbitrary ordering quantity $y\in{\cal Y}({\cal N})$ and an arbitrary ${\cal S}\subsetneq{\cal N}$ with ${\cal Y}({\cal S})=\mathbb{R}_+$. (Note that if ${\cal S}\subseteq {\cal N}_r$ for some $r$, $r=1,\ldots,R$, we then have ${\cal Y}({\cal S})=\{y^*({\cal S})\}$ and $v_{\max}(y,{\cal S})$ can be computed by simply solving a single optimization problem, $\displaystyle\min_{P\in{\cal P}(P_1,\ldots,P_R)}v_P(y,{\cal N})$.) We have:
\be
\label{eq:gvmax}
v_{\max}(y,{\cal S})=\max_{\gamma \in \R^+}~ \max_{P\in{{\cal P}(P_1,\ldots,P_R)}} \frac{(p-c)\gamma-p\,\mathbb{E}_P\left[\left(\gamma-\td{d}({\cal S})\right)^+\right]}
{(p-c)y-p\,\mathbb{E}_P\left[\left(y-\td{d}({\cal N})\right)^+\right]}.
\ee
For newsvendor games, it is reasonable to assume that demands follow discrete non-negative distributions, which could be constructed from historical sales or market analysis. More specifically, let each distribution $P_r$ of $\td{\mb{d}}_r$ be represented as a discrete non-negative distribution with $K_r$ values, $\mb{d}^k_r$ of probability $p^k_r$, for $k = 1,\ldots,K_r$, $r=1,\ldots,R$. Thus, each probability distribution $P$ in $\mathcal{P}(P_1,\ldots,P_R)$ is a discrete distribution with a support of $\displaystyle K =\prod_{r=1}^RK_r$ values $\mb{d}_k$ and each has an unknown probability of $q_k$, $k = 1,\ldots,K$. For $P$ to be consistent with $P_1,\ldots,P_R$, the following constraints on $\mb{q}$ must hold:
\be
\label{eq:qdef}
\left\{\displaystyle \mb{q} \geq 0,~~\sum_{k=1}^Kq_k=1,~~\sum_{k=1}^K\mathbb{I}\{\mb{d}_{k,r}=\mb{d}_r^l\}q_k = p_r^l,~~ r=1,\ldots,R,\,l=1,\ldots,K_r\right\}.
\ee
Given the one-to-one mapping between $P$ and $\mb{q}$, we abuse the notations and write both $v_P(y,{\cal S})$ and $v_{\vect{q}}(y,{\cal S})$ interchangeably when the context is clear for $\mb{q}$ to be the corresponding representation of $P$.

Problem \refs{eq:gvmax} can be reformulated as
\begin{equation}
\label{eq:vmax}
\begin{array}{rl}
v_{\max}(y,{\cal S})=\displaystyle \max_{\gamma, \vect{q}} & \displaystyle\frac{(p-c)\gamma - p\,\sum_{k=1}^K\left[\left(\gamma-d_k({\cal S})\right)^+\right]q_k}{(p-c)y - p\sum_{k=1}^K\left[\left(y-d_k({\cal N})\right)^+\right]q_k}\\
\st & \displaystyle\sum_{k=1}^K\mathbb{I}\{\mb{d}_{k,r}=\mb{d}_r^l\}q_k = p_r^l,\quad\forall\,r=1,\ldots,R,\,l=1,\ldots,K_r,\\
& \displaystyle\sum_{k=1}^Kq_k=1,\\
& \gamma,\mb{q}\geq 0.
\end{array}
\end{equation}
For each fixed $\gamma$, we can apply the standard method for transforming a linear fractional optimization problem into a linear program (see, for example, Cambini et al. \cite{cambini05}). To this end, let us introduce new decision variables
$$\displaystyle \theta = \frac{1}{(p-c)y - p\sum_{k=1}^K\left[\left(y-d_k({\cal N})\right)^+\right]q_k},$$
 and
$$\displaystyle \psi_k =  \frac{q_k}{(p-c)y - p\sum_{k=1}^K\left[\left(y-d_k({\cal N})\right)^+\right]q_k},\quad k=1,\ldots,K.$$
Under Assumption~\ref{as:pos}(ii), we have $\theta > 0$. The objective function then becomes $$\displaystyle (p-c)\gamma\cdot \theta - p\,\sum_{k=1}^K\left[\left(\gamma-d_k({\cal S})\right)^+\right]\psi_k,$$
and Problem~\refs{eq:vmax} can be reformulated as

\begin{equation}
\label{eq:mpriLP}
\begin{array}{rl}
v_{\max}(y,{\cal S})=\displaystyle \max_{\gamma, \theta, \vect{\psi}} & \displaystyle (p-c)\gamma\cdot \theta - p\,\sum_{k=1}^K\left[\left(\gamma-d_k({\cal S})\right)^+\right]\psi_k\\
\st & \displaystyle\sum_{k=1}^K\mathbb{I}\{\mb{d}_{k,r}=\mb{d}_r^l\}\psi_k - p_r^l\cdot\theta = 0,\quad\forall\,r=1,\ldots,R,\,l=1,\ldots,K_r,\\
& \displaystyle\sum_{k=1}^K \psi_k-\theta = 0,\\
& \displaystyle (p-c)y\cdot \theta - p\,\sum_{k=1}^K\left[\left(y-d_k({\cal N})\right)^+\right] \psi_k = 1,\\
& \gamma,\theta, \mb{\psi} \geq 0,
\end{array}
\end{equation}
where the first two constraints in \refs{eq:mpriLP} are derived directly from the first two constraints in \refs{eq:vmax} by multiplying both sides of those with $\theta$. The third constraint in \refs{eq:mpriLP} is derived by the definitions of $\theta$ and $\mb{\psi}$. Finally, we have replaced the constraint $\theta > 0$ by $\theta \geq 0$ without loss of generality since $\theta = 0$ is not a feasible solution (otherwise, $\mb{\psi}$ must be equal to zero from the second constraint and that violates the third constraint.)

Problem \refs{eq:mpriLP} is a bilinear optimization problem, which is generally not easy to solve. We show, however, in the following proposition that one of the distinct values of $\{d_1({\cal S}),\ldots,d_K({\cal S})\}$, which are known, is an optimal value of $\gamma$. Given a fixed $\gamma$, the resulting bilinear optimization problem is reduced to a linear program. It implies that we can solve Problem \refs{eq:mpriLP} by solving at most $K$ linear programs with $\gamma$ set to each and every distinct value of $\{d_1({\cal S}),\ldots,d_K({\cal S})\}$.

\begin{proposition}
\label{prop:linearequiv}
There exists an optimal solution $(\gamma^*,\theta^*,\mb{\psi}^*)$ of Problem \refs{eq:mpriLP} such that $$\gamma^*\in\{d_1({\cal S}),\ldots,d_K({\cal S})\}.$$
\end{proposition}
%The intuitive explanation of this result is that, for any fixed choice of $\mb{q}$, i.e., if the distribution $P\in{{\cal P}}$ is fixed, then an choice of the optimal order quantity is the $\displaystyle \frac{(p-c)}{p}$-quantile of $P(d({\cal S}))$. Given that $P(d({\cal S}))$ is discrete with the support set $\{d_1({\cal S}),\ldots,d_K({\cal S})\}$, that set also contains an optimal order quantity. We now provide a more formal proof of this result.

\begin{pf}
Given an arbitrary value of $\gamma$, Problem \refs{eq:mpriLP} is reduced to a linear program for $\theta$ and $\mb{\psi}$ over a fixed feasible set $\cal F$ defined by the set of constraints in \refs{eq:mpriLP}, i.e., $\gamma$ only affects the objective function. Under Assumption~\ref{as:pos}(ii), $\theta$ and $\mb{\psi}$ are non-negative and bounded, which means $\cal F$ is bounded and Problem \refs{eq:mpriLP} can be written as follows:
$$
v_{\max}(y,{\cal S}) = \max_{\gamma\geq 0}\left(\max_{s=1,\ldots,S}\left\{(p-c)\gamma\cdot \theta^s - p\,\sum_{k=1}^K\left[\left(\gamma-d_k({\cal S})\right)^+\right]\psi_k^s\right\}\right),
$$
where $\left\{(\theta^s,\mb{\psi}^s)\right\}_{s=1,\ldots,S}$ is the set of extreme points of $\cal F$. Equivalently, we have:
$$
v_{\max}(y,{\cal S}) = \max_{s=1,\ldots,S}\left\{\max_{\gamma\geq 0}\left((p-c)\gamma\cdot \theta^s - p\,\sum_{k=1}^K\left[\left(\gamma-d_k({\cal S})\right)^+\right]\psi_k^s\right)\right\}.
$$

For an arbitrary solution $(\theta^s,\mb{\psi}^s)$, $s=1,\ldots,S$, it is easy to show that function $\displaystyle f(\gamma;\theta^s,\mb{\psi}^s)=(p-c)\gamma\cdot \theta^s - p\,\sum_{k=1}^K\left[\left(\gamma-d_k({\cal S})\right)^+\right]\psi_k^s$ is a \emph{concave} piece-wise linear function with intersection points as distinct values of the set $\{d_1({\cal S}),\ldots,d_K({\cal S})\}$. Since $\theta^s>0$ as shown previously, $f(\cdot;\theta^s,\mb{\psi}^s)$ tends to $-\infty$ when $\gamma$ tends to $+\infty$. It means there is at least an optimal solution for the problem $\displaystyle\max_{\gamma\geq 0}f(\gamma;\theta^s,\mb{\psi}^s)$ which belongs to the set $\{d_1({\cal S}),\ldots,d_K({\cal S})\}$. Thus we have:
$$
v_{\max}(y,{\cal S}) = \max_{s=1,\ldots,S}\left\{\max_{l=1,\ldots,K}\left\{(p-c)d_l({\cal S})\cdot \theta^s - p\,\sum_{k=1}^K\left[\left(d_l({\cal S})-d_k({\cal S})\right)^+\right]\psi_k^s\right\}\right\},
$$
which shows that there exists an optimal solution $(\gamma^*,\theta^*,\mb{\psi}^*)$ of Problem \refs{eq:mpriLP} such that
 $$\gamma^*\in\{d_1({\cal S}),\ldots,d_K({\cal S})\}.$$

%Without loss of generality, let us first assume $\{d_1({\cal S}) < d_2({\cal S}) < \ldots < d_K({\cal S})\}$.
%For each fixed choice of $(\kappa, \mb{\psi})$, the objective function in Model~\ref{eq:mpriLP} is the summation of a linear function $(p-c)y({\cal S}) \kappa$ and $K$ piece-wise concave functions $- p\,\sum_{k=1}^K\left[\left(y({\cal S})-d_k({\cal S})\right)^+\right]\psi_k$ on $y({\cal S})$. For each $k =1,\ldots,K$, each of the piece-wise linear function $\left[\left(y({\cal S})-d_k({\cal S})\right)^+\right]\psi_k$ includes the first piece with a slope of zero and the second piece with a slope of $-p \psi_k$ and the joint is at $d_k({\cal S})$. Thus, the objective function is also a piece-wise concave function on $y({\cal S})$ with pieces joining each other at $d_k({\cal S}), k=1,\ldots, K$.

%The constraint set in Model~\ref{eq:mpriLP} forms a polyhedron on $(\kappa, \mb{\psi})$ with a finite number of extreme points and extreme rays, each of which results in a piece-wise concave function on $y({\cal S})$ as explained above. Therefore, the maximum of these piece-wise concave functions is another piece-wise function on $y({\cal S})$ which could be neither convex nor concave. However, since the piece-wise concave functions all share the same set of joints at $\{d_1({\cal S}),\ldots,d_K({\cal S})\}$, one of these joints attain the maximum for the objective function.
\end{pf}

Proposition \ref{prop:linearequiv} shows us how to compute $v_{\max}(y,{\cal S})$ by solving at most $K$ linear programs. We can use this approach to compute $v_{\max}(y^*_{wc}({\cal N}),{\cal S})$ as inputs of the linear program \refs{eq:rnleastcore}, which is then solved to check the existence of core solutions of the robust newsvendor game $({\cal N},{\cal Y},{\cal V}({\cal P}(P_1,\ldots,P_r)))$. If the core is empty, we need to consider the general problem \refs{eq:rleastcore} applying to the robust newsvendor game, whose optimal solutions can be considered as its least core solutions. Let us consider the following problem, which is similar to \refs{eq:rnleastcore}, for an arbitrary $y \in {\cal Y}({\cal N})$:
\be
\label{eq:ly}
\ba{rl}
\displaystyle \sigma(y)=\min_{\bx,\eps} & \eps\\
\st & \displaystyle\sum_{i\in{\cal S}}x_i\geq v_{\max}(y,{\cal S})-\eps,\quad{\cal S}\subsetneq{\cal N},\\
& \displaystyle\sum_{i\in{\cal N}}x_i=1.
\ea
\ee
\begin{comment}
\begin{equation}
\label{eq:ly}
\begin{array}{rl}
\displaystyle \Upsilon(y) = \min_{\mb{x},\epsilon} & \epsilon\\
s.t. & \mb{x}({\cal S}) + \epsilon \geq v_{max}(y,{\cal S}), \forall {\cal S} \subsetneq {\cal N},\\
& \mb{e}^T \mb{x} = 1.
\end{array}
\end{equation}
\end{comment}
Clearly, $\displaystyle s({\cal N},{\cal Y},{\cal V}({\cal P}(P_1,\ldots,P_r)))=\min_{y \in {\cal Y}({\cal N})} \sigma(y)$, which is a reformulation of the least core problem \refs{eq:rleastcore}.  We will show that $\sigma(y)$ is a convex function in the following proposition.
\begin{comment}
We can utilize the convexity property of $\Upsilon(y)$ in the numerical computation of a $y^*$ with the smallest least core value. Specifically, we can search for the (one dimensional) $y$ that is a local optimal solution and conclude that this is also the global optimal solution.
\end{comment}
\begin{proposition}
\label{prop:ly_vmax_properties}
The following statements hold:
\begin{itemize}
\item[(a)] For each coalition ${\cal S}\subsetneq {\cal N}$, $\displaystyle v_{max}(y,{\cal S})$ is a convex function of $y$ on ${\cal Y}({\cal N})$.
\item[(b)] $\sigma(y)$ is a convex function of $y$ on ${\cal Y}({\cal N})$.
\end{itemize}
\end{proposition}

In order to prove the proposition, we need the following lemma.

\begin{lemma}
\label{lemma:vmax_properties}
The inverse function $\displaystyle \frac{1}{v_P(y,{\cal N})}$ is a convex function of $y$ on ${\cal Y}({\cal N})$ for all $P\in{\cal P}(P_1,\ldots,P_R)$.
\end{lemma}
%The intuitive explanation of this result is that, for any fixed choice of $\mb{q}$, i.e., if the distribution $P\in{{\cal P}}$ is fixed, then an choice of the optimal order quantity is the $\displaystyle \frac{(p-c)}{p}$-quantile of $P(d({\cal S}))$. Given that $P(d({\cal S}))$ is discrete with the support set $\{d_1({\cal S}),\ldots,d_K({\cal S})\}$, that set also contains an optimal order quantity. We now provide a more formal proof of this result.

\begin{pf}
Let $\mb{q}$ be the corresponding probability vector of a joint distribution $P\in {\cal P}(P_1,\ldots,P_R)$. We have:
$$
v_P(y,{\cal N})\equiv v_{\vect{q}}(y,{\cal N})=(p-c)y - p\sum_{k=1}^K\left[\left(y-d_k({\cal N})\right)^+\right]q_k,
$$
which is a piecewise linear concave function of $y$ with at most $(K+1)$ linear pieces positioning within the intervals induced by the sorted sequence of $\{d_1({\cal N}),d_2({\cal N})\ldots d_K({\cal N})\}$. We can therefore rewrite $\displaystyle v_P(y,{\cal N})=\min_{k=1,\ldots,K+1}\{a_ky+b_k\}$ where $a_k,b_k$ are appropriate linear coefficients that can be derived from $p$, $c$, $\mb{d}({\cal N})$ and $\mb{q}$. Since $v_P(y,{\cal N})>0$ for all $y\in{\cal Y}({\cal N})$, we have: $a_ky+b_k>0$ for all $k=1,\ldots,K+1$, and $y\in{\cal Y}({\cal N})$. We then have:
$$
\frac{1}{v_P(y,{\cal N})}=\max_{k=1,\ldots,K+1}\frac{1}{a_ky+b_k},
$$
which is the maximum of convex inverse linear functions on ${\cal Y}({\cal N})$ and hence is also a convex function in ${\cal Y}({\cal N})$.
\begin{comment}
Let $d^{(1)}({\cal N})\leq d^{(2)}({\cal N})\leq \ldots \leq d^{(K)}({\cal N})$ be the corresponding sorted sequence of $\{d_1({\cal N}),d_2({\cal N})\ldots d_K({\cal N})\}$. We also denote $q^{(j)}$ as the corresponding probability of the realization of the joint demand at $d^{(j)}({\cal N})$ for each $j = 1,\ldots,K$.

Let us denote $a_k = (\rho-c - \rho \sum_{j=1}^k q^{(j)})$ and $b_k = \rho \sum_{j=1}^k q^{(j)} d^{(k)}({\cal N})$ for $k= 0,\ldots,K+1$. We have
$$
\displaystyle v_P(y,{\cal N}) =
\begin{cases}
\displaystyle a_0 y + b_0, &\mbox{if } y \leq d^{(1)}({\cal N}) \\
\displaystyle a_k y+b_k, &\mbox{if } d^{(k)}({\cal N}) \leq y \leq d^{(k+1)}({\cal N})\\
\displaystyle a_K y+b_K, &\mbox{if } y \geq d^{(K)}({\cal N}).
\end{cases}
$$
Let us define $\displaystyle g_k(y) = \frac{v_P(\gamma,{\cal S})}{a_k y+b_k},~k=0,\ldots,K+1,$ which are inverse linear functions and are convex. Then $\displaystyle \frac{v_P(\gamma,{\cal S})}{v_P(y,{\cal N})}$ composes of these $K+1$ convex pieces. In addition, by the definition of $(a_k,b_k),~k=0,\ldots,K$, we can show that
$$ \displaystyle v_P(y,{\cal N}) = \max_{j=0,\ldots,K} g_j(y),$$
which is the maximum of convex functions and hence is also a convex function.
\end{comment}
\end{pf}

We are now ready to prove the Proposition~\ref{prop:ly_vmax_properties}.

\begin{pf}
(a) If ${\cal S}\subseteq{\cal N}_r$ for some $r$, $r=1,\ldots,R$, we have: ${\cal Y}({\cal S})=\{y^*({\cal S})\}$ and
$$
\ba{rl}
v_{\max}(y,{\cal S})&=\displaystyle\max_{P\in{{\cal P}(P_1,\ldots,P_R)}} \frac{\bar{v}({\cal S})}
{(p-c)y-p\,\mathbb{E}_P\left[\left(y-\td{d}({\cal N})\right)^+\right]}\\
& =\displaystyle\frac{\bar{v}({\cal S})}{\min_{\vect{q}\in{\cal Q}}\left\{(p-c)y - p\sum_{k=1}^K\left[\left(y-d_k({\cal N})\right)^+\right]q_k\right\}}\equiv \displaystyle\frac{\bar{v}({\cal S})}{\min_{\vect{q}\in{\cal Q}} v_{\vect{q}}(y,{\cal N})},
\ea
$$
where ${\cal Q}={\cal Q}(P_1,\ldots,P_R)$ is the feasible set of the probability vector $\mb{q}$ as described in \refs{eq:qdef}. The second equality follows from the fact that $\bar{v}({\cal S})\geq 0$ and $v_P(y,{\cal N}) > 0$. Since ${\cal Q}$ is a bounded polytope; there exists an optimal solution $\mb{q}^*\in{\cal Q}^*$, where ${\cal Q}^*$ is the set of extreme points of ${\cal Q}$. Thus we have:
$$
v_{\max}(y,{\cal S})=\max_{\vect{q}\in{\cal Q}^*}\frac{\bar{v}({\cal S})}{v_{\vect{q}}(y,{\cal N})}.
$$
Since $\bar{v}({\cal S})\geq 0$ and since $\displaystyle \frac{1}{v_{\vect{q}}(y,{\cal N})}$ is convex for each $\mb{q}$ according to Lemma \ref{lemma:vmax_properties}, we have: $v_{\max}(y,{\cal S})$ is convex.

Now consider an arbitrary ${\cal S}\subsetneq{\cal N}$ with ${\cal Y}({\cal S})=\mathbb{R}_+$. We have,
$$
v_{\max}(y,{\cal S})=\max_{P\in{{\cal P}(P_1,\ldots,P_R)}} \frac{\displaystyle\max_{\gamma \in {\cal Y}({\cal S})}v_P(\gamma,{\cal S})}
{v_P(y,{\cal N})}.
$$
We have: $v_P(y,{\cal N})>0$ for all $P\in{\cal P}(P_1,\ldots,P_R)$ and $y\in{\cal Y}({\cal N})$. In addition, ${\cal Y}({\cal S})=\mathbb{R}_+$, thus there always exists $\gamma\in{\cal Y}({\cal S})$ small enough such that $v_P(\gamma,{\cal S})\geq 0$ for any $P$. We then have: $v_{\max}(y,{\cal S})\geq 0$ for all $y\in{\cal Y}({\cal N})$. We can rewrite the formulation of $v_{\max}(y,{\cal S})$ as follows:
$$
v_{\max}(y,{\cal S})=\max_{\gamma \in \R^+}~ \max_{\vect{q}\in{{\cal Q}}} \frac{(p-c)\gamma - p\sum_{k=1}^K\left[\left(\gamma-d_k({\cal S})\right)^+\right]q_k}
{(p-c)y - p\sum_{k=1}^K\left[\left(y-d_k({\cal N})\right)^+\right]q_k}.
$$

Proposition~\ref{prop:linearequiv} shows that we can restrict the domain of $\gamma$ to the discrete set $\{d_1({\cal S}),\ldots,d_K({\cal S})\}$, that is,
$$
v_{\max}(y,{\cal S})=\max_{\gamma \in \{d_1({\cal S}),\ldots,d_K({\cal S})\}}~ \max_{\vect{q}\in{{\cal Q}}} \frac{(p-c)\gamma - p\sum_{k=1}^K\left[\left(\gamma-d_k({\cal S})\right)^+\right]q_k}
{(p-c)y - p\sum_{k=1}^K\left[\left(y-d_k({\cal N})\right)^+\right]q_k}.
$$
\begin{comment}
Thus, Problem~\refs{eq:mpriLP} can be reformulated as

\begin{equation}
\label{eq:vmax2}
\begin{array}{rl}
v_{\max}(y,{\cal S})=\displaystyle \max_{\gamma \in \{d_1({\cal S}),\ldots,d_K({\cal S})\}} ~~ \max_{\vect{q} \in Q} & h(y,{\cal S},\gamma,\mb{q}),
\end{array}
\end{equation}
where $\displaystyle h(y,{\cal S},\gamma,\mb{q}) = \displaystyle\frac{(p-c)\gamma - p\,\sum_{k=1}^K\left[\left(\gamma-d_k({\cal S})\right)^+\right]q_k}{(p-c)y - p\sum_{k=1}^K\left[\left(y-d_k({\cal N})\right)^+\right]q_k}$.
\end{comment}
For a fixed $\gamma$, the inner problem is a linear fractional optimization problem over the bounded polyhedron $\cal Q$ and hence there exists an optimal solution $\mb{q}^*$ given that $v_P(y,{\cal N})>0$ for all $y\in{\cal Y}({\cal N})$ and $P\in{\cal P}(P_1,\ldots,P_R)$. Let us consider the level sets ${\cal L}_{\alpha}$ of the linear fractional objective function, which are hyperplanes. For the optimal objective value $\alpha^*$, we have: ${\cal L}_{\alpha^*}\cap{\cal Q}\neq\emptyset$. We claim that ${\cal L}_{\alpha^*}\cap{\cal Q}^*\neq\emptyset$, where ${\cal Q}^*$ is the set of extreme points of $\cal Q$. Since $\alpha^*$ is the optimal objective value, $\cal Q$ belongs to a half-space defined by ${\cal L}_{\alpha^*}$. Suppose, on contradiction, that ${\cal L}_{\alpha^*}\cap{\cal Q}^*=\emptyset$. Due to convexity and since the entire ${\cal Q}^*$ belongs to the same half-space defined by ${\cal L}_{\alpha^*}$, we have: ${\cal L}_{\alpha^*}\cap{\cal Q}=\emptyset$ (contradiction).

With ${\cal L}_{\alpha^*}\cap{\cal Q}^*\neq\emptyset$, we can now compute $v_{\max}(y,{\cal S})$ as follows:
$$
v_{\max}(y,{\cal S})=\max_{\gamma \in \{d_1({\cal S}),\ldots,d_K({\cal S})\}}~ \max_{\vect{q}\in{{\cal Q}^*}} \frac{(p-c)\gamma - p\sum_{k=1}^K\left[\left(\gamma-d_k({\cal S})\right)^+\right]q_k}
{(p-c)y - p\sum_{k=1}^K\left[\left(y-d_k({\cal N})\right)^+\right]q_k}.
$$
Since $v_{\max}(y,{\cal S})\geq 0$, we can focus on the set ${\cal H} \in \{d_1({\cal S}),\ldots,d_K({\cal S})\} \times {\cal Q}^*$ of $(\gamma,\mb{q})$ such that $(p-c)\gamma - p\sum_{k=1}^K\left[\left(\gamma-d_k({\cal S})\right)^+\right]q_k\geq 0$ when computing $v_{\max}(y,{\cal S})$, that is,
$$
v_{\max}(y,{\cal S})=\max_{(\gamma,\vect{q}) \in {\cal H}} \frac{(p-c)\gamma - p\sum_{k=1}^K\left[\left(\gamma-d_k({\cal S})\right)^+\right]q_k}
{(p-c)y - p\sum_{k=1}^K\left[\left(y-d_k({\cal N})\right)^+\right]q_k}.
$$
Applying Lemma \refs{lemma:vmax_properties}, clearly, $v_{\max}(y,{\cal S})$ is the maximum of convex functions, which means it is also a convex function.
\begin{comment}
Let $\mb{q}^{(1)},\ldots,\mb{q}^{(M)}$ be the corresponding extreme points. Since level sets of the linear fractional objective function are hyperplanes, the optimal level set that passes through $\mb{q}^* \in Q$ must also pass through at least an extreme point of $Q$. This is because otherwise all these extreme points must belongs to the same half-plane produced by the optimal level set hyperplane, but not on that hyperplane, which means the entire polyhedron $Q$ does not intersect with the hyperplane and this contradicts with $\mb{q}^* \in Q$.

Thus, the optimality of Problem~\refs{eq:vmax2} can be attained at one of the extreme points $\mb{q}^{(1)},\ldots,\mb{q}^{(M)}$.  Problem~\refs{eq:vmax2} can then be reformulated as
\begin{equation}
\label{eq:vmax3}
\begin{array}{rl}
v_{\max}(y,{\cal S})=\displaystyle \max_{\gamma \in \{d_1({\cal S}),\ldots,d_K({\cal S})\}} ~~ \max_{\vect{q} \in \{\mb{q}^{(1)},\ldots,\mb{q}^{(M)}\}} & h(y,{\cal S},\gamma,\mb{q}),
\end{array}
\end{equation}
which is the maximum of a number of functions $h(y,{\cal S},d_l({\cal S}),\mb{q}^{(m)}),~l=1,\ldots,K,~m=1,\ldots,M$. Since each of these function is convex over $y$ by Lemma~\ref{lemma:vmax_properties}, we have $v_{\max}(y,{\cal S})$ is also a convex function of $y$.
\end{comment}

(b) To prove the convexity of $\sigma(y)$, we show that, for any $\{y_1, y_2, y_3\} \in {\cal Y}({\cal N})$ such that there exists $\alpha \in [0,1]$ with $y_2 = \alpha y_1+ (1-\alpha) y_3$, then $\sigma(y_2) \leq \alpha \sigma(y_1) + (1-\alpha) \sigma(y_3)$.

Let $(\mb{x}_1,\epsilon_1)$ and $(\mb{x}_3,\epsilon_3)$ be the optimal solutions of ~\refs{eq:ly} when $y = y_1$ and $ y=y_3$, respectively. Let us define $(\mb{x}_2,\epsilon_2) = \alpha (\mb{x}_1,\epsilon_1)+(1-\alpha)(\mb{x}_3,\epsilon_3)$. It is easy to verify that $\mb{e}^T \mb{x}_2 = 1$. In addition, for all $ {\cal S} \subsetneq {\cal N}$, we have
\begin{eqnarray}
\mb{x}_2({\cal S}) + \epsilon_2 &=& \alpha (\mb{x}_1({\cal S})+\epsilon_1)+(1-\alpha)(\mb{x}_3({\cal S})+\epsilon_3) \label{eq:ly1}\\
 &\geq& \alpha v_{max}(y_1,{\cal S})+(1-\alpha) v_{max}(y_3,{\cal S}) \label{eq:ly2}\\
 &\geq& v_{max}( \alpha y_1+ (1-\alpha) y_3,{\cal S}) \label{eq:ly3}\\
 &=& v_{max}( y_2,{\cal S}) \label{eq:ly4},
\end{eqnarray}
where \refs{eq:ly1} comes directly from the construction of $(\mb{x}_2,\epsilon_2)$; \refs{eq:ly2} comes from the feasibility of $(\mb{x}_1,\epsilon_1)$ and $(\mb{x}_3,\epsilon_3)$; \refs{eq:ly3} comes from the convexity of $v_{max}(y,{\cal S})$ as shown in part (a). Finally, \refs{eq:ly4} comes directly from the definition of $y_2$.

This shows that $(\mb{x}_2,\epsilon_2)$ is a feasible solution of ~\refs{eq:ly} when $y = y_2$. Therefore,
$$\sigma(y_2) \leq \epsilon_2= \alpha \epsilon_1 + (1-\alpha) \epsilon_3= \alpha \sigma(y_1) + (1-\alpha) \sigma(y_3),$$
i.e., $\sigma(y)$ is a convex function.
\end{pf}

%Without loss of generality, let us first assume $\{d_1({\cal S}) < d_2({\cal S}) < \ldots < d_K({\cal S})\}$.
%For each fixed choice of $(\kappa, \mb{\psi})$, the objective function in Model~\ref{eq:mpriLP} is the summation of a linear function $(p-c)y({\cal S}) \kappa$ and $K$ piece-wise concave functions $- p\,\sum_{k=1}^K\left[\left(y({\cal S})-d_k({\cal S})\right)^+\right]\psi_k$ on $y({\cal S})$. For each $k =1,\ldots,K$, each of the piece-wise linear function $\left[\left(y({\cal S})-d_k({\cal S})\right)^+\right]\psi_k$ includes the first piece with a slope of zero and the second piece with a slope of $-p \psi_k$ and the joint is at $d_k({\cal S})$. Thus, the objective function is also a piece-wise concave function on $y({\cal S})$ with pieces joining each other at $d_k({\cal S}), k=1,\ldots, K$.

%The constraint set in Model~\ref{eq:mpriLP} forms a polyhedron on $(\kappa, \mb{\psi})$ with a finite number of extreme points and extreme rays, each of which results in a piece-wise concave function on $y({\cal S})$ as explained above. Therefore, the maximum of these piece-wise concave functions is another piece-wise function on $y({\cal S})$ which could be neither convex nor concave. However, since the piece-wise concave functions all share the same set of joints at $\{d_1({\cal S}),\ldots,d_K({\cal S})\}$, one of these joints attain the maximum for the objective function.

Proposition \ref{prop:ly_vmax_properties} shows that the least core problem \refs{eq:rleastcore} for our robust newsvendor game is a convex optimization problem in terms of $y$ and we could apply simple one-dimensional search algorithms to find the optimal solution. The next section provides some numerical results on the properties and computation of the core (and least core) solutions of robust newsvendor games.
\subsection{Numerical Results}
%\todo{Numerical examples, compare with the joint independent demand distribution}
%In this section, we compare the performance of the robust allocation scheme (called ROBUST) with its deterministic counterpart (called INDEPT and to be defined in details later). To this end, we assume the following experimental setting.
We consider the following experimental setting. We are given a set of retailers ${\cal N}$ and a partition ${\cal N}_1,\ldots,{\cal N}_R$. In addition, for each $r= 1,\ldots, R$, we are given the discrete historical joint demand distribution $P_r$ for the subset of retailers ${\cal N}_r$ but not the joint demand distribution of all retailers. Discussions in Section~\ref{robust_core_computation} allow us to compute a robust core (or least core) solution by solving \refs{eq:rnleastcore} (and \refs{eq:rleastcore} if necessary) under the framework of robust newsvendor games, which consists of the allocation scheme $\bz_{rob}$ and an order quantity $y_{rob}$ for the grand coalition. %as well as the order quantity $y_{rob}({\cal S})$ for each group of retailers ${\cal S} \subsetneq {\cal N}$ in the case the group wants to form the coalition itself.

In order to evaluate the performance of the robust solution $(y_{rob},\bz_{rob})$, we are going to compare it with the solution derived from the deterministic newsvendor game under the assumption that all multivariate marginal distributions $P_r$, $r=1,\ldots,R$, are independent of each other, that is, $\displaystyle \mathbb{P}(\tilde{\mb{d}} = (\mb{d}_1^{l_1},\ldots,\mb{d}_r^{l_R})) = \prod_{r=1}^R \mathbb{P}(\tilde{\mb{d}}_r = \mb{d}_r^{l_r})=\prod_{r=1}^R p_r^{l_r}$. This could be considered as a common assumption on the joint distribution given its marginal distributions. Clearly, the resulting joint distribution $P_I$ belongs to $\mathcal{P}(P_1,\ldots,P_R)$. Given this distribution $P_I$, we can compute the allocation scheme $\bz_{det}=\bx/\bar{v}({\cal N})$, where $\bx$ is a core solution of the deterministic newsvendor game with respect to $P_I$. In addition, the optimal order quantity $y_{det}$ for the grand coalition in this deterministic newsvendor game is used to form the solution $(y_{det},\bz_{det})$, which will be compared with the robust solution $(y_{rob},\bz_{rob})$. %$y_{det}({\cal N})$ for the grand coalition as well as the order quantity $y_{det}({\cal S})$ for each group of retailers ${\cal S} \subsetneq {\cal N}$ in the case the group wants to form the coalition itself.

We shall compare the performance of these two solutions with respect to joint distributions which belong to $\mathcal{P}(P_1,\ldots,P_R)$. Given a distribution $P\in \mathcal{P}(P_1,\ldots,P_R)$, we compute the maximum normalized dissatisfaction or worst normalized excess value for each solution. For $(y_{rob},\bz_{rob})$, the excess value is computed as
\be
\label{eq:excess}
\eps_{rob}^P=\max_{{\cal S}\subsetneq{\cal N}}\left\{\left(\frac{\displaystyle\max_{\gamma\in{\cal Y}({\cal S})}v_P(\gamma,{\cal S})}{v_P(y_{rob},{\cal N})}-z_{rob}({\cal S})\right)^+\right\}.
\ee
The excess value $\eps_{det}^P$ can be defined in the same fashion for $(y_{det},\bz_{det})$. We follow the stress test approach proposed by Dupa\v cov\' a \cite{dupacova06} with the contaminated distributions $P_{\lambda} = \lambda P_I + (1-\lambda)P_{ext}$ for $\lambda\in[0,1]$, where $P_{ext}$ are extremal distributions, i.e., those distributions which are likely to be the ones with which $v_{\max}(y,{\cal S})$ are computed. Similar approach has been discussed in Bertsimas et al. \cite{bertsimas10} to test the quality of some stochastic optimization solutions.

We now consider a numerical example with $n=10$ and $R=2$, with the sizes of subsets, $\card{{\cal N}_1} = 4$ and $\card{{\cal N}_2} = 6$, respectively. We construct multivariate marginal distributions $P_1$ and $P_2$ from a randomly generated discrete joint distribution with the support set of each individual retailer's demand set to $[1,10]$. Other parameters include $p=1.5$ and $c= 1$. All the numerical results are tested on a PC with 2.67 gigahertz CPU, 12 gigabyte RAM, and a 64-bit Windows 7 operating system. We use MATLAB 8.0 for coding and IBM CPLEX Studio Academic version 12.5 for solving LPs problems under default settings. In order to compute a robust core (or least core) solution of this game, we would need to compute $v_{\max}(y,{\cal S})$ for each of $2^n=1024$ coalitions. This is accomplished by solving a number of linear programs as presented in Proposition~\ref{prop:linearequiv}. On average, the total time it took to compute a single value $v_{\max}(y,{\cal S})$ is approximately $12$ seconds under this setting. Problem \refs{eq:rnleastcore} can then be directly solved whereas \refs{eq:rleastcore} is solved with a simple one-dimensional search algorithm whose main subroutine depends on the solution of \refs{eq:rnleastcore} for different values of $y$.

\begin{figure}[htp]
 \begin{center}
\includegraphics[width=0.7\textwidth]{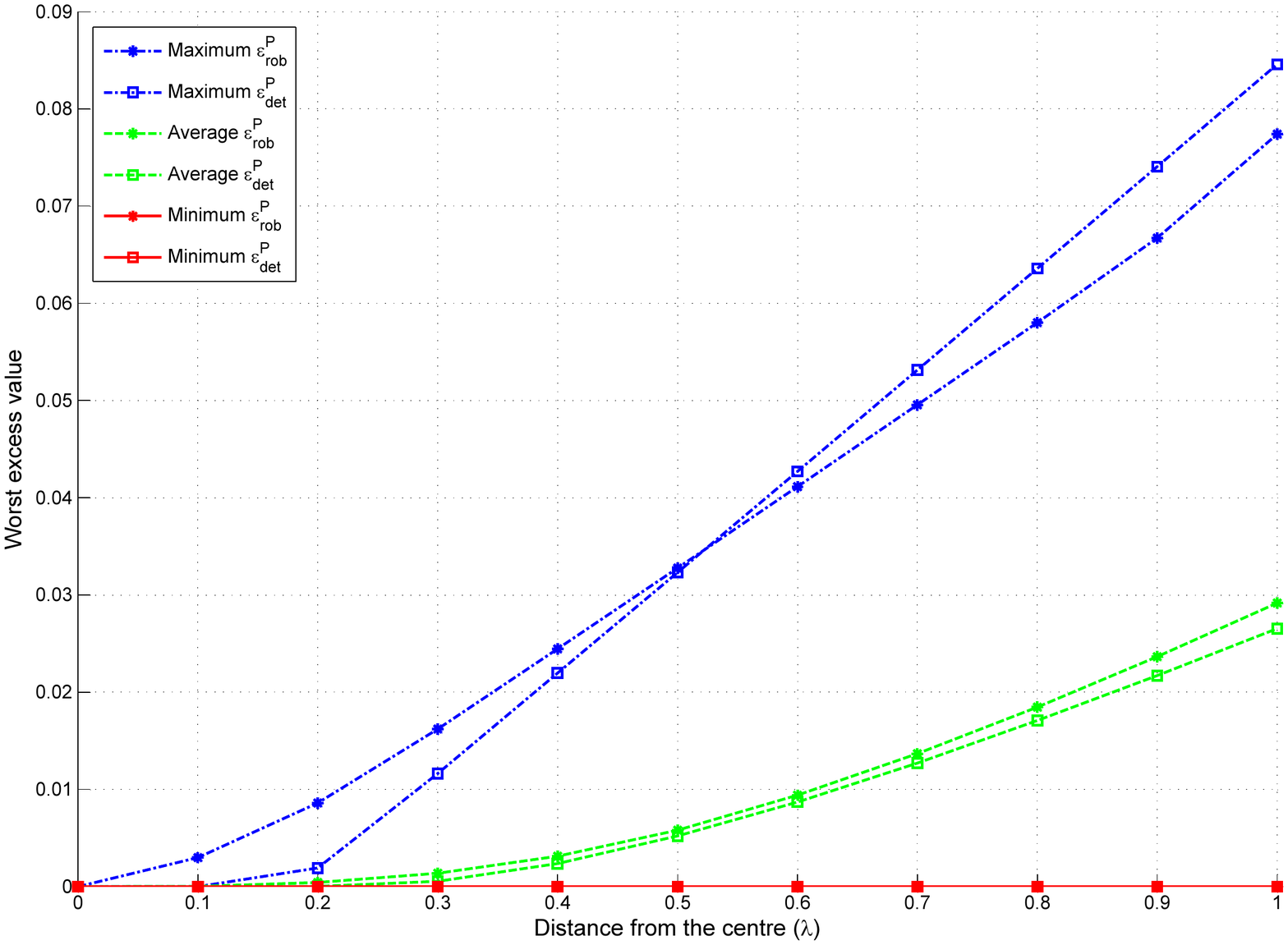}
\end{center}
\caption{Comparison between $\eps_{rob}^P$ and $\eps_{det}^P$ for contaminated distributions with different $\lambda$.}
\label{fig:robust_vs_deterministic}
\end{figure}

In this numerical example, for a given value of $\lambda$, we simply generate $100$ random extremal distributions $P_{ext}$ by solving the linear program $\{ \max \mb{c}^T \mb{q} ~:~ \mb{q} \in {\cal Q}\}$ with random cost vectors $\mb{c}$. We also include extremal distributions produced while calculating $v_{max}$. The excess values $\eps_{rob}^P$ and $\eps_{det}^P$ are computed for all contaminated distributions $P_\lambda$. Figure \ref{fig:robust_vs_deterministic} provides comparisons on three statistics of $\eps_{rob}^P$ and $\eps_{det}^P$ for each fixed value of $\lambda$: the maximum, the minimum, and the average.
The maximum values of both $\eps_{rob}^P$ and $\eps_{det}^P$ increase when $\lambda$ increases. For $\lambda>0.5$, the robust solution $(y_{rob},\bz_{rob})$ yields smaller excess values in the worst case as compared to those of $(y_{det},\bz_{det})$ for contaminated distributions. It shows that the robust solution hedges against the worst case as expected even though on average, the solution $(y_{det},\bz_{det})$ is slightly better in terms of worst excess values. It is worth noting that in the best case, both solutions are core solutions with no dissatisfaction even for the case of $\lambda=1$.

\begin{figure}[htp]
 \begin{center}
\includegraphics[width=0.7\textwidth]{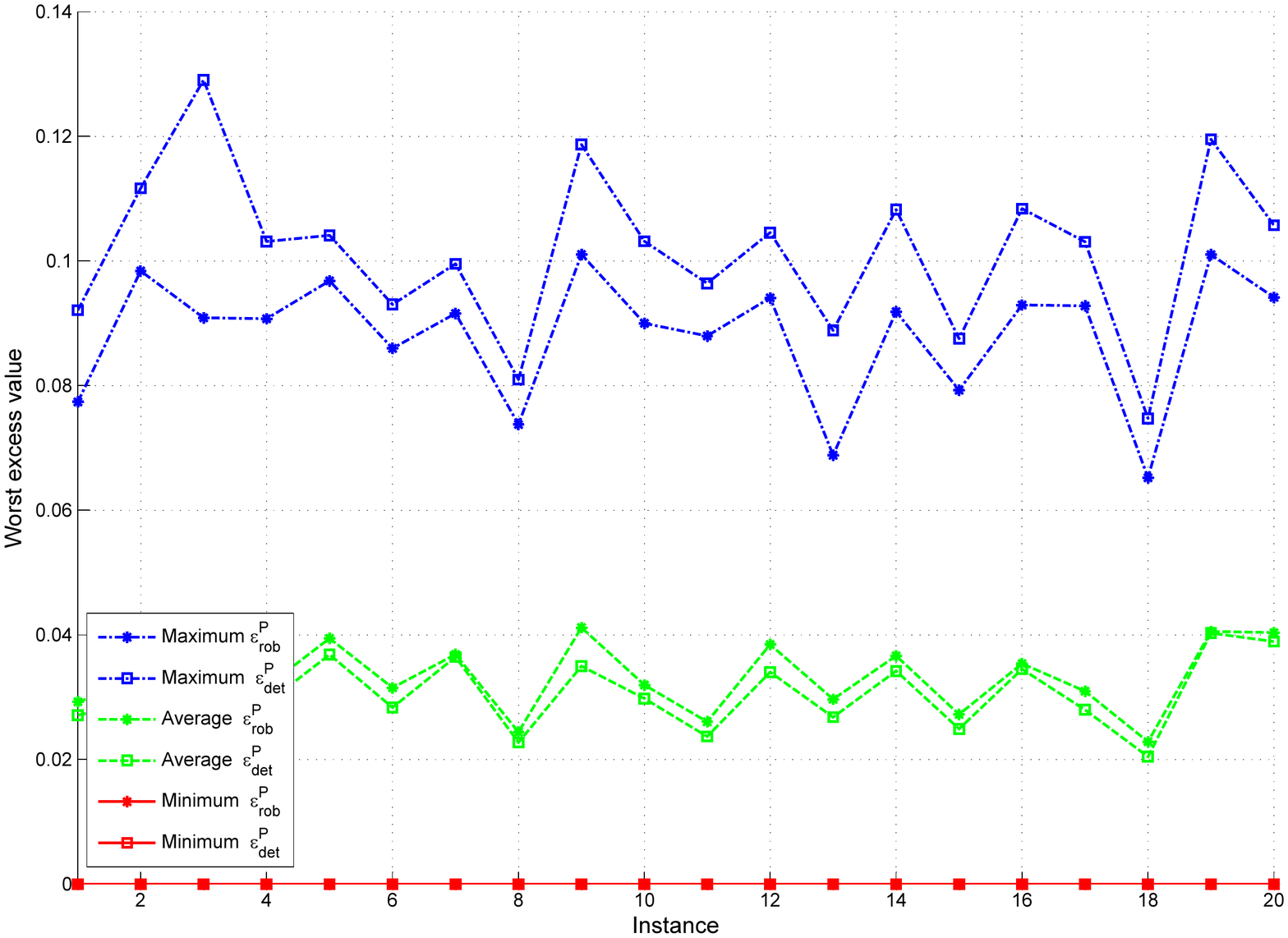}
\end{center}
\caption{Comparison between $\eps_{rob}^P$ and $\eps_{det}^P$ for different instances with $\lambda=1$.}
\label{fig:robust_vs_deterministic1}
\end{figure}

\begin{figure}[htp]
 \begin{center}
\includegraphics[width=0.7\textwidth]{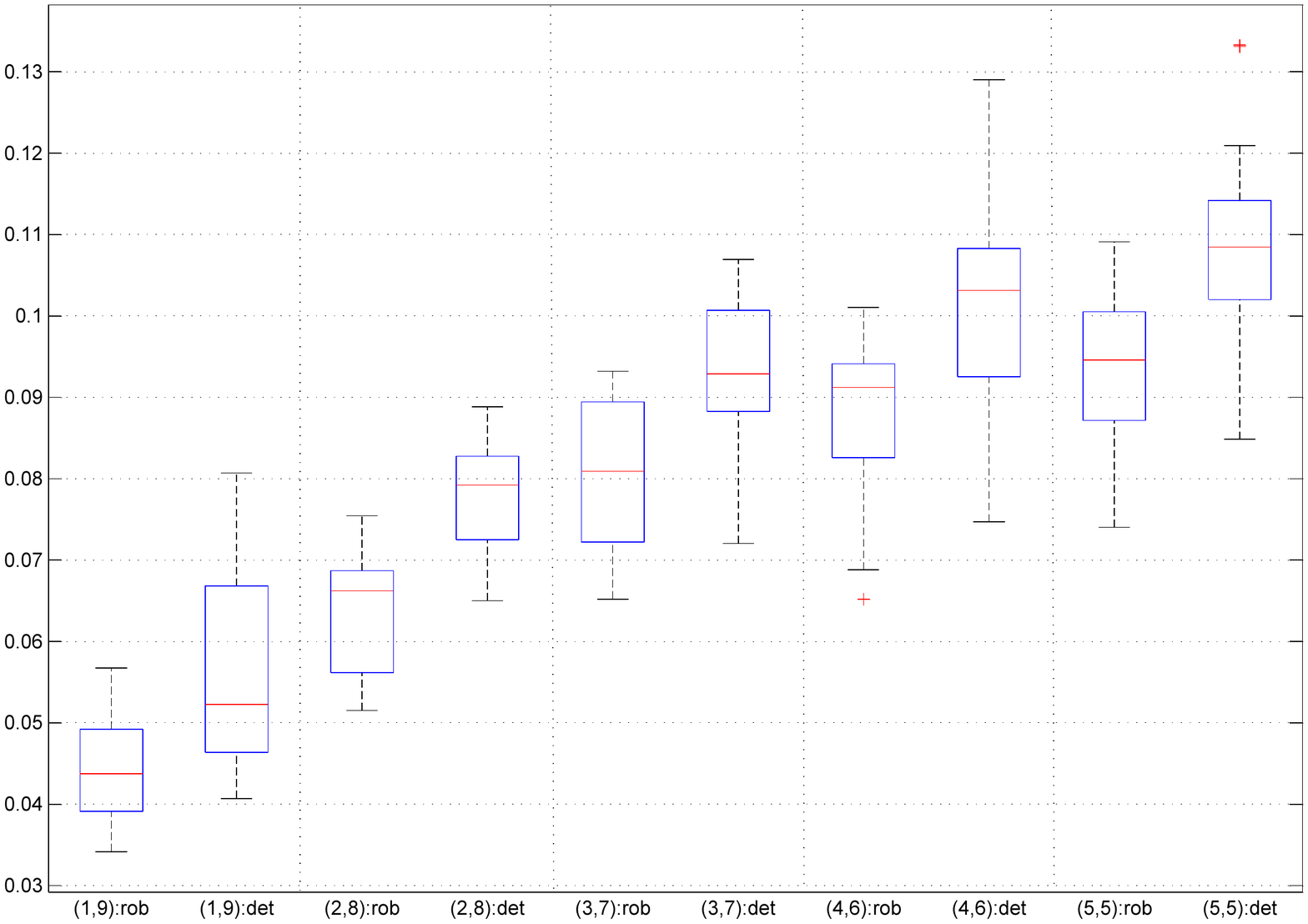}
\end{center}
\caption{Comparison between maximum $\eps_{rob}^P$ and $\eps_{det}^P$ for different subset structures with $\lambda=1$.}
\label{fig:robust_vs_deterministic2}
\end{figure}

We run the experiment again for $M=20$ different instances. Figure \ref{fig:robust_vs_deterministic1} shows the statistics of $\eps_{rob}^P$ and $\eps_{det}^P$ when $\lambda=1$ for all of these instances. The results again show that the robust solution $(y_{rob},\bz_{rob})$ consistently outperforms $(y_{det},\bz_{det})$ for all these instances in the worst case and is slightly worse on average.
%on larger games with $n=10$ retailers. We consider a setting with $R=2$ and with ${\cal N}_1 = \{1,2,3,4,5\}$, ${\cal N}_2 = \{6,7,8,9,10\}$. In order to find the robust core/least core of this game, we need to solve the bilinear Problem~\refs{eq:mpriLP} for each of $(2^n=1024)$ coalitions. This is accomplished by solving at most $K$ linear programs as presented in Proposition~\ref{prop:linearequiv}. On average, the total time it took to solve Problem~\refs{eq:mpriLP} is around 12 seconds. All the numerical results are tested on a personal computer with 2.67GHz CPU, 12GB RAM, and the Windows 7, 64-bit operating system. We use MATLAB for coding and use IBM CPLEX Studio Academic version 12.5 for solving LPs problems under default settings. Comparison between ROBUST and INDEPT are shown in Figure~\ref{fig:robust_vs_deterministic4} in Appendix B where a similar trend as observed for the case of $n=3$ is presented.
Finally, we run the experiment again for four different settings with respect to sizes of the two subsets, $(1,9)$, $(2,8)$, $(3,7)$, and $(5,5)$, in addition to the original setting of $(4,6)$, with $M=20$ instances for each setting. Figure \ref{fig:robust_vs_deterministic2} shows the box plots for the maximum values of $\eps_{rob}^P$ and $\eps_{det}^P$ when $\lambda=1$. The results again show that the robust solution $(y_{rob},\bz_{rob})$ outperforms $(y_{det},\bz_{det})$ in the worst case.
\section{Conclusion}
In this paper, we develop a framework for newsvendor games with ambiguity in demand distributions, which we call robust cooperative games. We discuss solution concepts of robust cooperative games and study them in the context of newsvendor games with ambiguity in demand distributions when only marginal distributions are known. Some numerical results are provided, which show the robust core solutions hedge against the worst cases as expected. It is possible to develop other frameworks for cooperative games with uncertain characteristic functions by using different payoff distribution schemes and preference relations, which could be applied to other applications.

%Our numerical results show attractive properties of the robust payoff distribution in comparison to non-robust counterparts. Specifically, we show that, using the suggested robust payoff distribution, coalitions are more stable to perturbations of the distribution function.  We demonstrate how the stability value and the degree of uncertainty affect the robust core existence through numerical simulations. We also show the algorithm developed can handle reasonably large stochastic newsvendor games.
%\bibliographystyle{plain}
\bibliographystyle{plainnat}
\bibliography{SNG}

\section*{Appendix: Derivation for Example \ref{ex:nocore}}\label{app:example1}
Let us consider the robust newsvendor game introduced in Example \ref{ex:nocore}. According to Lemma \ref{lem:optquant}, the worst-case optimal ordering quantity $y_{wc}^*({\cal N})$ is computed as follows:
$$
y_{wc}^*({\cal N}) = y^*(\{1,2\})+ y^*(\{3\}) = D+\frac{D(p-c)}{p}=\frac{D(2p-c)}{p}.
$$
Similarly,
$$
v_{wc}({\cal N}) = \bar{v}(\{1,2\}) + \bar{v}(\{3\}) = D(p-c)+\frac{D(p-c)^2}{2p}=\frac{D(p-c)(3p-c)}{2p}.
$$
Given $P\in{\cal P}(P_1,P_2)$, we have:
\begin{eqnarray*}
v_P(y^*_{wc}({\cal N}),{\cal N}) &=& (p-c)(y^*(\{1,2\})+y^*(\{3\})) - p\mathbb{E}_{P}(y^*(\{1,2\})+y^*(\{3\}) - D - \td{d}_3)^+\\
&=& (p-c)D+(p-c)y^*(\{3\}) - p\mathbb{E}_{P_2}(y^*(\{3\}) -\td{d}_3))^+\\
&=& \bar{v}(\{1,2\}) + \bar{v}(\{3\})\\
&=& v_{wc}({\cal N}).
\end{eqnarray*}
Thus we have:
\begin{eqnarray*}
v_{max}(y^*_{wc}({\cal N}),{\cal S}) &=& \displaystyle \max_{P\in{\cal P}(P_1,P_2)} ~ \max_{y\in{\cal Y}({\cal S})} ~\frac{v_{P}(y,{\cal S})}{v_P(y^*_{wc}({\cal N}),{\cal N})}\\
&=& \displaystyle \max_{P\in{\cal P}(P_1,P_2)} ~ \max_{y\in{\cal Y}({\cal S})}~ \frac{v_{P}(y,{\cal S})}{v_{wc}({\cal N})}.
\end{eqnarray*}
Combining with \refs{eq:thm:rcore1}, we can compute $v_{\max}(y^*_{wc}({\cal N}),\{1,2\})$ as follows:
$$
v_{max}(y^*_{wc}({\cal N}),\{1,2\}) = \frac{\bar{v}(\{1,2\})}{v_{wc}(N)} = \frac{2p}{3p-c}.
$$
Now consider the sub-coalition $\{1,3\}$. Let $P$ be the distribution characterized by the marginal distribution of $\td{d}_1$ and the relationship $\td{d}_2=\td{d}_3=D-\td{d}_1$. Clearly, $P\in{\cal P}(P_1,P_2)$. In addition, let $y=D$, which is a feasible ordering quantity, we have:
\begin{eqnarray*}
v_{max}(y^*_{wc}({\cal N}),\{1,3\}) &\geq& \frac{v_{P}(y,\{1,3\})}{v_{wc}({\cal N})}\\
&=& \frac{D(p-c) - p\mathbb{E}_{P}(D - D)^+}{v_{wc}({\cal N})} \\
&=& \frac{2p}{3p-c}.
\end{eqnarray*}
Given the symmetry between $\td{d}_1$ and $\td{d}_2$, we can show that $v_{max}(y^*_{wc}({\cal N}),\{2,3\})=v_{max}(y^*_{wc}({\cal N}),\{1,3\})$. Thus, we obtain the required inequalities,
$$
v_{\max}(y^*_{wc}({\cal N}),\{1,2\})=\frac{2p}{3p-c}\leq v_{\max}(y^*_{wc}({\cal N}),\{2,3\})=v_{\max}(y^*_{wc}({\cal N}),\{1,3\}),
$$
for Example \ref{ex:nocore}.
\begin{comment}
\section*{Appendix B: Comparison between robust and deterministic solution}\label{app:other_numerical_results}
\begin{figure}[htp]
 \begin{center}
 %\includegraphics[scale=1]{equal_weight1.jpg}
% \includegraphics[scale=0.45]{SNG_n3R1R2.eps}
\end{center}
\caption{Comparison between robust solutions and deterministic solutions (with $n=3,~R=2$ and with ${\cal N}_1 = \{1\}$, ${\cal N}_2 = \{2,3\}$.)}
\label{fig:robust_vs_deterministic2}
\end{figure}

\begin{figure}[htp]
 \begin{center}
 %\includegraphics[scale=1]{equal_weight1.jpg}
% \includegraphics[scale=0.45]{SNG_n3R1R1R1.eps}
\end{center}
\caption{Comparison between robust solutions and deterministic solutions (with $n=3,~R=3$ and with ${\cal N}_1 = \{1\}$, ${\cal N}_2 = \{3\}$, ${\cal N}_3 = \{3\}$.)}
\label{fig:robust_vs_deterministic3}
\end{figure}

\begin{figure}[htp]
 \begin{center}
 %\includegraphics[scale=1]{equal_weight1.jpg}
% \includegraphics[scale=0.45]{SNG_n10R5R5.eps}
\end{center}
\caption{Comparison between robust solutions and deterministic solutions (with $n=10,~R=2$ and with ${\cal N}_1 = \{1,2,3,4,5\}$, ${\cal N}_2 = \{6,7,8,9,10\}$.)}
\label{fig:robust_vs_deterministic4}
\end{figure}
\end{comment}
\end{document}